\documentclass{amsart}
\usepackage{graphics}
\usepackage{amssymb}
\usepackage[all]{xy}
\CompileMatrices
\newtheorem*{introthm}{Theorem}

\newtheorem*{introcor}{Corollary}
\newtheorem{thm}{Theorem}[section]
\newtheorem{lemma}[thm]{Lemma}
\newtheorem{prop}[thm]{Proposition}
\newtheorem{coro}[thm]{Corollary}
\theoremstyle{definition}
\newtheorem*{introdefi}{Definition}
\newtheorem*{introexam}{Example}
\newtheorem{defi}[thm]{Definition}
\newtheorem{exam}[thm]{Example}
\newtheorem{rem}[thm]{Remark}
\newtheorem{constr}[thm]{Construction}
\numberwithin{equation}{thm}

\def\Rho{{\rm R}}

\def\div{{\rm div}}
\def\quot{/\!\!/}

\def\rq#1{\widehat{#1}}
\def\t#1{\widetilde{#1}}
\def\b#1{\overline{#1}}

\def\C{{\rm C}}

\def\CC{{\mathbb C}}
\def\KK{{\mathbb K}}

\def\ZZ{{\mathbb Z}}

\def\QQ{{\mathbb Q}}
\def\PP{{\mathbb P}}
\def\CDiv{\operatorname{CDiv}}
\def\WDiv{\operatorname{WDiv}}
\def\PDiv{\operatorname{PDiv}}

\def\Cl{\operatorname{Cl}}
\def\ClDiv{\operatorname{Cl}}
\def\Pic{\operatorname{Pic}}
\def\Hom{{\rm Hom}}

\def\codim{{\rm codim}}

\def\Spec{{\rm Spec}}

\def\cone{{\rm cone}}
\def\faces{{\rm faces}}

\def\lin{{\rm lin}}
\def\cov{{\rm cov}}

\def\rank{\operatorname{rank}}

\def\topto#1{\stackrel{{\scriptscriptstyle #1}}{\longrightarrow}}
\begin{document}
\title[Bunches of cones]
      {Bunches of cones in the divisor \\ 
       class group --- A new combinatorial \\
       language for toric varieties}
\author[F. Berchtold]{Florian Berchtold} 
\address{Fachbereich Mathematik und Statistik, Universit\"at Konstanz,
78457 Konstanz, Germany}
\email{Florian.Berchtold@uni-konstanz.de}
\author[J.~Hausen]{J\"urgen Hausen} 
\address{Mathematisches Forschungsinstitut Oberwolfach, Lorenzenhof,
77709 Oberwolfach-Walke, Germany}
\email{hausen@mfo.de}
\begin{abstract}
As an alternative to the description of a toric variety
by a fan in the lattice of one parameter subgroups, 
we present a new language in terms of what we call
bunches --- these are certain collections of cones in the
vector space of rational divisor classes. 
The correspondence between these bunches and fans is 
based on classical Gale duality. 
The new combinatorial language allows a much more natural 
description of geometric phenomena around divisors of 
toric varieties than the usual method by fans does. 
For example, the numerically effective cone and the ample cone
of a toric variety can be read off immediately from its bunch. 
Moreover, the language of bunches appears to be useful for
classification problems.
\end{abstract}

\maketitle

\section*{Introduction}

The most important feature of a toric variety $X$ is that it is
completely described by a fan $\Delta$ in the lattice $N$ of 
one parameter subgroups of the big torus $T_{X} \subset X$. 
Applying a linear Gale transformation to the set of
primitive generators $v_{\varrho}$ of the rays 
$\varrho \in \Delta^{(1)}$ of $\Delta$, gives a new
vector configuration in a rational vector space $K_{\QQ}$.
This opens an alternative combinatorial approach to the
toric variety $X$: The vector space $K_{\QQ}$ is isomorphic
to the rational divisor class group of $X$, and one can 
shift combinatorial information between the spaces $N_{\QQ}$ 
and $K_{\QQ}$.

In toric geometry, this principle has been used to study 
the projective case, compare~\cite{OP}: Roughly speaking, 
if we consider all fans $\Sigma$ in $N$ having their rays among 
$\Delta^{(1)}$, then the (quasi-)projective $\Sigma$ correspond 
to so called Gelfand-Kapranov-Zelevinski cones in $K_{\QQ}$.
These cones subdivide the cone generated by the Gale transform
of the vector configuration 
$\{v_{\varrho}; \; \varrho \in \Delta^{(1)}\}$,
and the birational geometry of the associated toric varieties is
reflected by the position of their Gelfand-Kapranov-Zelevinski 
cones.

If one leaves the (quasi-)projective setting, then there are
generalizations of Gelfand-Kapranov-Zelevinski cones,
compare e.g.~\cite{Hu}; but so far there seems to be no concept which is 
simple enough to serve for practical purposes in toric geometry.  
Our aim is to fill this gap and to propose a natural
combinatorial language which also works in the non 
quasiprojective case.
The combinatorial data are certain collections
 --- which we call {\em bunches} --- of overlapping cones 
in the vector space of rational divisor classes.
As we shall see, this approach gives very natural 
descriptions of geometric phenomena connected with divisors.

\goodbreak

In order to give a first impression of the language of bunches, let 
us present it here for in a special case, namely for toric
varieties with free class group. 
Consider a sequence $(w_{1}, \dots, w_{n})$ of not necessarily
pairwise distinct points, called {\em weight vectors}, in a lattice
$K \cong \ZZ^{k}$.
By a {\em weight cone\/} in $K$ we then mean a (convex polyhedral)
cone $\tau \subset K_{\QQ}$ generated by some of the $w_{i}$.

\begin{introdefi}
A {\em  (free standard) bunch\/} in $K$ is a nonempty
collection $\Theta$ of weight cones in $K$ with the following 
properties: 
\begin{enumerate}
\item a weight cone $\sigma$ in $K$ belongs to $\Theta$ if and 
  only if 
   $\emptyset \ne \sigma^{\circ} \cap \tau^{\circ} \ne \tau^{\circ}$
   holds for all $\tau \in \Theta$ with $\tau \ne \sigma$.
\item for each $i$, the set $\{w_{j}; \; j \ne i\}$ generates $K$ 
  as a lattice, and there is a $\tau \in \Theta$ such that 
  $\tau^{\circ} \subset \cone(w_{j}; \; j \ne i)^{\circ}$
  holds for the relative interiors.
\end{enumerate}
\end{introdefi}

How to construct a toric variety from such a bunch $\Theta$? The first
step is to {\em unpack\/} the combinatorial information encoded in
$\Theta$. For this, let $E := \ZZ^{n}$, and let $Q \colon E \to K$ 
denote the linear surjection sending the canonical base vector 
$e_{i}$ to $w_{i}$.
Denote by $\gamma \subset E_{\QQ}$ the positive orthant. We define the 
{\em covering collection\/} of $\Theta$ as
$$ \cov(\Theta) := \{\gamma_{0} \preceq \gamma; \; 
                   \gamma_{0} \text{ minimal with } 
                   Q(\gamma_{0}) \supset \tau \text{ for some } 
                   \tau \in \Theta\}. $$

The next step is to {\em dualize\/} the information contained in
$\cov(\Theta)$. This is done by a procedure close to a
linear Gale transformation, for the classical setup see e.g.~\cite{Gr} 
and~\cite{MM}: Consider the exact sequence
arising from $Q \colon E \to K$, and the corresponding dualized exact
sequence:
$$ \xymatrix{
    0 \ar[r] &
    {M} \ar[r] &
    {E} \ar[r]^{Q} & 
    K \ar[r] &  
    0,  }$$
$$ \xymatrix{
    0 \ar[r] &
    {L} \ar[r] &
    {F} \ar[r]^{P} & 
    {N} \ar[r] &  
    0 . }$$

Note that $P \colon F \to N$ is not the dual homomorphism of $Q \colon E \to
K$. Let $\delta := \gamma^{\vee} \subset F_{\QQ}$ be the dual cone of
$\gamma \subset E_{\QQ}$. 
The crucial observation then is that we obtain a
fan $\Delta(\Theta)$ in the lattice $N$ 
having as its maximal cones the images
$$
P(\gamma_{0}^{\perp} \cap \delta), \quad \gamma_{0} \in
\cov(\Theta). 
$$

\begin{introdefi}
The {\em toric variety associated to the bunch\/} $\Theta$
is $X_{\Theta} := X_{\Delta(\Theta)}$.
\end{introdefi}

The toric variety $X_{\Theta}$ is {\em nondegenerate\/} in the sense
that it has no torus factors. Moreover, $X_{\Theta}$ is 
{\em 2-complete},
that means if $X_{\Theta} \subset X$ is an open toric embedding 
such that the
complement $X \setminus X_{\Theta}$ is of codimension at least 2,
then $X = X_{\Theta}$.  

\begin{introexam}
Consider the sequence $(1,2,3)$ of weight vectors in $K := \ZZ$
and the bunch $\Theta := \{\QQ_{\ge 0}\}$ in $\ZZ$.
Then we have $E = \ZZ^{3}$, and the associated linear map 
$Q \colon \ZZ^{3} \to \ZZ$ sends $e_{i}$ to $i$.
The covering collection $\cov(\Theta)$ consists of
the following three faces of $\gamma = \QQ_{\ge 0}^{3}$:
$$ 
\gamma_{i} := \QQ_{\ge 0}e_{i},
\qquad i = 1,2,3.
$$
If we identify the dual space $F = E^{*}$ with $\ZZ^{3}$,
then the cone $\delta = \gamma^{\vee}$ is again $\QQ_{\ge 0}^{3}$. 
Moreover, we may identify $N$ with $\ZZ^{2}$ and thus realize 
the map $P \colon F \to N$ via the matrix
$$
\left[
\begin{array}{lll}
-2 & 1 & 0 \\
-3 & 0 & 1 
\end{array}
\right]
$$
Each $\gamma_{i}^{\perp} \cap \delta$ equals
$\cone(e_{j},e_{l})$, where $j \ne i$ and $l \ne i$.
Hence the images $P(\gamma_{i}^{\perp} \cap \delta)$,
which are the maximal cones of the fan $\Delta(\Theta)$,
are given in terms of the canonical base vectors 
$e_{1},e_{2} \in \ZZ^{2}$ as
$$ 
\cone(e_{1},e_{2}),
\qquad
\cone(-2e_{1}-3e_{2},e_{2}),
\qquad
\cone(-2e_{1}-3e_{2},e_{1}).
$$
Consequently, the toric variety $X_{\Theta}$ associated 
to the bunch $\Theta$ equals the weighted projective space 
$\PP_{1,2,3}$. 
\end{introexam}

Introducing a suitable notion of a morphism, we can extend the assignment
$\Theta \mapsto X_{\Theta}$ from bunches to 2-complete nondegenerate toric
varieties to a contravariant functor. In fact, we even obtain a weak 
antiequivalence, see Theorem~\ref{standard2maximal}:

\begin{introthm}
The functor $\Theta \mapsto X_{\Theta}$ induces a bijection on the
level of isomorphism classes of bunches and nondegenerate 2-complete
toric varieties. In particular, every complete toric variety
arises from a bunch. 
\end{introthm}

In order to read off geometric properties
of $X_{\Theta}$ directly from the bunch $\Theta$,
one has to translate the respective fan-theoretical formulations
via the above Gale transformation into the language of bunches.
This gives for example: 
\begin{itemize}
\item $X_{\Theta}$ is $\QQ$-factorial if and only if every 
  cone $\tau \in \Theta$ is of full dimension, see 
  Proposition~\ref{Qfactorialchar}.
\item $X_{\Theta}$ is smooth if and only if for 
  every $\gamma_{0}$ in $\cov(\Theta)$, the image 
  $Q(\lin(\gamma_{0}) \cap E)$ equals $K$,
  see Proposition~\ref{smoothchar}.
\item The orbits $X_{\Theta}$ have fixed points in their closures if
  and only if all cones $Q(\gamma_{0})$, where $\gamma_{0} \in \cov(\Theta)$, 
  are simplicial, see Proposition~\ref{existenoffixedpointschar}. 
\end{itemize}

As mentioned, the power of the language of bunches lies in the 
description of geometric phenomena around divisors, because $K_{\QQ}$ 
turns out to be the rational divisor class group 
of $X_{\Theta}$.
For example, we obtain very simple descriptions 
for the classes of rational Cartier divisors, 
the cone $\C^{\rm sa}(X_{\Theta})$ of
semiample classes
and the cone $\C^{\rm a}(X_{\Theta})$ of ample classes,
see Theorem~\ref{QCartier}:

\begin{introthm}
For the toric variety $X_{\Theta}$ arising from a bunch $\Theta$,
there are canonical isomorphisms:
$$
\Pic_{\QQ}(X_{\Theta}) \cong \bigcap_{\tau \in \Theta} \lin(\tau),
 \qquad
\C^{\rm sa}(X_{\Theta}) \cong \bigcap_{\tau \in \Theta} \tau,
 \qquad
\C^{\rm a}(X_{\Theta}) \cong \bigcap_{\tau \in \Theta} \tau^{\circ}.
$$
\end{introthm}

Note that the last isomorphism gives a quasiprojectivity criterion 
in the spirit of~\cite{Ew} and~\cite{Sh}, see Corollary~\ref{Ewaldcrit}.
Moreover, we can derive from the above Theorem 
a simple Fano criterion, see Corollary~\ref{fanocrit}.
Finally, we get back Reid's Toric Cone Theorem, see~\cite{Re},
even with a new description of the Mori Cone, 
see Corollary~\ref{MoriCone}:

\begin{introcor}
Suppose that $X_{\Theta}$ is complete and simplicial.
Then the cone of numerically effective 
curve classes in $H_{2}(X,\QQ)$ is given by
$$ 
{\rm NE}(X_{\Theta}) \cong  \sum_{\tau \in \Theta} \tau^{\vee}.
$$
In particular, this cone is convex and polyhedral. Moreover,
$X_{\Theta}$ is projective if and only if ${\rm NE}(X_{\Theta})$
is strictly convex.
\end{introcor}

Bunches can also be used for classification problems. 
For example, once the machinery is established, Kleinschmidt's
classification~\cite{Kl} becomes very simple
and can even be slightly improved, 
see Proposition~\ref{Kleinschmidt};
below we represent a sequence of weight vectors 
as a set of vectors $w$, each of which carries a
multiplicity $\mu(w)$: 

\begin{introthm}
The smooth 2-complete toric varieties $X$ with $\Cl(X) \cong \ZZ^{2}$
correspond to bunches $\Theta = \{\cone(w_{1},w_{2})\}$ given by
\begin{itemize}
\item weight vectors $w_{1} := (1,0)$, and 
$w_{i} := (b_{i},1)$ with $0 = b_{n} < b_{n-1} < \cdots < b_{2}$, 
\item multiplicities $\mu_{i} := \mu(w_{i})$ with $\mu_{1} > 1$,
  $\mu_{n} > 0$ and $\mu_{2} + \cdots + \mu_{n} > 1$. 
\end{itemize}

\medskip

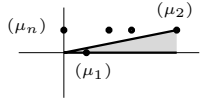
\begin{center}
  \begin{picture}(0,0)%
\includegraphics{smoothclass.pstex}%
\end{picture}%
\setlength{\unitlength}{1243sp}%
\begingroup\makeatletter\ifx\SetFigFont\undefined%
\gdef\SetFigFont#1#2#3#4#5{%
  \reset@font\fontsize{#1}{#2pt}%
  \fontfamily{#3}\fontseries{#4}\fontshape{#5}%
  \selectfont}%
\fi\endgroup%
\begin{picture}(3847,1937)(1126,-2333)
\put(2476,-1861){\makebox(0,0)[lb]{\smash{\SetFigFont{6}{7.2}{\familydefault}{\mddefault}{\updefault}
$(\mu_{1})$
}}}
\put(1126,-961){\makebox(0,0)[lb]{\smash{\SetFigFont{6}{7.2}{\familydefault}{\mddefault}{\updefault}
$(\mu_{n})$
}}}
\put(4101,-636){\makebox(0,0)[lb]{\smash{\SetFigFont{6}{7.2}{\familydefault}{\mddefault}{\updefault}
$(\mu_{2})$
}}}
\end{picture}

\end{center}
Moreover, the toric variety $X$ defined by such a bunch $\Theta$
is always projective, and it is Fano if and only if we have
$$ b_{2}(\mu_{3} + \cdots + \mu_{n}) < \mu_{1} + b_{2}\mu_{3} +
\cdots + b_{n-1}\mu_{n-1}. $$
\end{introthm}

In general the functor from bunches to
toric varieties is neither injective nor surjective on morphisms,
see Examples~\ref{freegroup} and~\ref{notfaithfull}. 
But if we restrict to $\QQ$-factorial toric varieties, then the
language of bunches provides also a tool for the study of toric
morphisms, see Theorem~\ref{equivcat}: 

\begin{introthm}
There is an equivalence from the category of simple bunches to
the category of full $\QQ$-factorial toric varieties.
\end{introthm}

\tableofcontents

We would like to thank the referee for his careful 
reading and for various helpful comments and corrections.

\section{Fans and toric varieties}\label{section1}

In this section, we recall some basic facts on the correspondence 
between fans and toric varieties, and thereby fix our notation 
used later. 
For details, we refer to the books of Oda~\cite{Od} and 
Fulton~\cite{Fu}.
We begin with introducing the necessary terminology from convex 
geometry.

By a lattice we mean a free finitely generated $\ZZ$-module. The
associated rational vector space of a lattice $N$ is $N_{\QQ} := \QQ
\otimes_{\ZZ} N$. If $P \colon F \to N$ is a
lattice homomorphism, then we denote the induced linear map 
$F_{\QQ} \to N_{\QQ}$ of rational vector spaces again by $P$.

By a cone in a lattice $N$ we always mean a polyhedral (not
necessarily strictly) convex cone in
the associated rational vector space $N_{\QQ}$. 
Let $N$ be a lattice, and let $M
:= \Hom(N,\ZZ)$ denote the dual lattice of $N$. The orthogonal space
and the dual cone of a cone $\sigma$ in $N$ are 
$$ 
\sigma^{\perp} := \{u \in M_{\QQ}; \; u\vert_{\sigma} = 0\},
\qquad
\sigma^{\vee} := \{u \in M_{\QQ}; \; u\vert_{\sigma} \ge 0\}.$$

The relative interior of a cone $\sigma$ 
is denoted by $\sigma^{\circ}$.
If $\sigma_{0}$ is a face of
$\sigma$, then we write $\sigma_{0} \preceq \sigma$. 
The dimension of $\sigma$ is the dimension of 
the linear space $\lin(\sigma)$ generated by $\sigma$. 
The set of the $k$-dimensional faces of $\sigma$ is denoted by 
$\sigma^{(k)}$, 
and the one-dimensional faces of a strictly convex
cone are called rays.

The primitive generators of a strictly convex
cone $\sigma$ in a lattice $N$ are the
primitive lattice vectors of its rays. 
A strictly convex cone in $N$ is called 
simplicial if its primitive generators are linearly independent,
and it is called regular if its primitive generators can
be complemented to a lattice basis of $N$.

\begin{defi}\label{fandefi}
\begin{enumerate}
\item A {\em fan\/} in a lattice $N$ is a finite collection $\Delta$
  of strictly convex cones in $N$ such that for each $\sigma \in
  \Delta$ also all $\sigma_{0} \preceq \sigma$ belong to $\Delta$ and
  for any two $\sigma_{i} \in \Delta$ we have 
  $\sigma_{1} \cap \sigma_{2}  \preceq \sigma_{i}$. 
\item A {\em map of fans\/} $\Delta_{i}$ in lattices $N_{i}$ is a
  lattice homomorphism $F \colon N_{1} \to N_{2}$ such that for every
  $\sigma_{1} \in \Delta_{1}$ there is a $\sigma_{2} \in \Delta_{2}$
  containing $F(\sigma_{1})$.
\end{enumerate}
\end{defi}

Recall that the compatibility condition 
$\sigma_{1} \cap \sigma_{2}  \preceq \sigma_{i}$
in the above definition is equivalent to the existence of a 
{\em separating linear form\/} for the cones $\sigma_{1}$
and $\sigma_{2}$, i.e., a linear form $u$ on $N$
such that 
$$
u_{\vert \sigma_{1}} \ge 0,
\qquad
u_{\vert \sigma_{2}} \le 0,
\qquad
u^{\perp} \cap \sigma_{i} = \sigma_{1} \cap \sigma_{2}.
$$

If we replace in Definition~\ref{fandefi} ``strictly convex'' with
``convex'', we obtain the
category of {\em quasifans\/}. For a fan $\Delta$ in $N$, we denote by
$\vert \Delta \vert$ its support, that is the union of all its
cones. Moreover, $\Delta^{\max}$ is the set of maximal cones of
$\Delta$, and $\Delta^{(k)}$ is the set of all $k$-dimensional cones of
$\Delta$.  

In the sequel, we shall often make use of a well known 
universal lifting construction, which makes the set of primitive 
generators of the rays of a given fan into a lattice basis, 
compare for example~\cite{Co}:

\begin{constr}\label{coxconstr}
Let $\Delta$ be a fan in a lattice $N$, 
and let $\Rho := \Delta^{(1)}$. 
Let $C \colon \ZZ^{\Rho} \to N$ be the map sending the
canonical base vector $e_{\varrho}$ to the primitive generator
$v_{\varrho} \in \varrho$. For $\sigma \in \Delta^{\max}$, set
$$ \t{\sigma} := \cone(e_{\varrho}; \; \varrho \in \sigma^{(1)}). $$  
Then the cones $\t{\sigma}$, where $\sigma \in \Delta^{\max}$, are the
maximal cones of a fan $\t{\Delta}$ consisting of faces of the
positive orthant in $\QQ^{\Rho}$. Moreover, $C \colon \ZZ^{\Rho} \to N$
is a map of the fans $\t{\Delta}$ and $\Delta$.
\end{constr}

Now we turn to toric varieties. Throughout the entire paper, we work
over an algebraically closed field $\KK$ of characteristic zero, and
the word point refers to a closed point.

\begin{defi}
\begin{enumerate}
\item A {\em toric variety\/} is a normal variety $X$ containing an
  algebraic torus $T_{X}$ as an open subset such that the group
  structure of $T_{X}$ extends to a regular action on $X$.
\item A {\em toric morphism\/} is a regular map $X \to Y$ of toric
  varieties that restricts to a group homomorphism $T_{X} \to T_{Y}$.
\end{enumerate}
\end{defi}

The correspondence between fans and toric varieties 
is obtained as follows: Let $\Delta$ be a fan in a lattice $N$,
and let $M := \Hom(N,\ZZ)$ be the dual lattice of $N$. 
For every cone $\sigma \in \Delta$ one
defines an affine toric variety:
$$ X_{\sigma} := \Spec(\KK[\sigma^{\vee} \cap M]). $$

For any two such $X_{\sigma_{i}}$, one has canonical open 
embeddings of $X_{\sigma_{1} \cap \sigma_{2}}$ into $X_{\sigma_{i}}$. 
Patching together all $X_{\sigma}$ along these open
embeddings gives a toric variety $X_{\Delta}$.
The assignment $\Delta \mapsto X_{\Delta}$ is functorial;
it is even a (covariant) equivalence of categories.

In the sequel, we shall frequently restrict our investigations
to toric varieties that behave reasonably. 
For that purpose, we consider the following geometric properties:

\begin{defi}\label{nondegfulldef}
\begin{enumerate}
\item A toric variety $X$ is called {\em nondegenerate\/} if it admits
  no toric decomposition $X \cong X' \times \KK^{*}$.
\item We call a toric variety $X$ {\em 2-complete\/} if it does not
  admit a toric open embeding $X \subset X'$ with $X' \setminus X$
  nonempty of codimension at least two.
\item We call a toric variety $X$ {\em full\/} if it is
  2-complete and every $T_{X}$-orbit has a fixed point
  in its closure. 
\end{enumerate}
\end{defi}

The notion of 2-completeness already occurs in~\cite{BeHa}.
It generalizes completeness in the sense that a toric variety
is complete if and only if it is ``1-complete''.  
For example, the affine space $\KK^{n}$ is 2-complete,
whereas for a toric variety $X$ of dimension at least two
and a fixed point $x \in X$, the variety $X \setminus \{x\}$ 
is not 2-complete.

In terms of fans, the properties introduced in 
Definition~\ref{nondegfulldef} are characterized as
follows:

\begin{rem}
Let $X$ be the toric variety arising from a fan $\Delta$ in a lattice
$N$. 
\begin{enumerate}
\item $X$ is nondegenerate if and only if the support $\vert \Delta
  \vert$ generates $N_{\QQ}$ as a vector space.
\item $X$ is 2-complete if and only if the fan $\Delta$ cannot
  be enlarged without adding new rays.
\item $X$ is full if and only if $\Delta$ is as in (ii) and every
  maximal cone of $\Delta$ is of full dimension.
\end{enumerate}
\end{rem}

\section{The category of bunches}~\label{section2}

In this section, we introduce the language of bunches.
Intuitively, one should think of a bunch as a
collection of pairwise overlapping lattice cones,
which satisfies certain irredundancy and maximality
properties.

The precise definition of the category of bunches is performed
in three steps.
The first one is to introduce the category of projected
cones:

\begin{defi}\label{weightedlatticedef}
\begin{enumerate}
\item A {\em projected cone\/} is a pair $(E \topto{Q} K, \gamma)$,
  where $Q \colon E \to K$ is an epimorphism of lattices and 
  $\gamma \subset E_{\QQ}$ is a simplicial cone of
  full dimension. 
\item A {\em morphism of projected cones\/} 
  $(E_{i} \topto{Q_{i}} K_{i}, \gamma_{i})$
  is a homomorphism $\Phi \colon E_{1} \to 
  E_{2}$ such that $\Phi(\gamma_{1}) \subset \gamma_{2}$ holds and 
  there is a commutative diagram
$$ \xymatrix{
E_{1} \ar[d]_{Q_{1}} \ar[r]^{\Phi} 
& 
E_{2} \ar[d]^{Q_{2}} \\
K_{1} \ar[r]_{\b{\Phi}} & K_{2} 
}$$ 
\end{enumerate}
\end{defi}

In the second step, we give the definition of bunches.
Such a bunch will live in a projected cone
$(E \topto{Q} K, \gamma)$.
By a {\em projected face\/} in $K$ we mean the image
$Q(\gamma_{0})$ of a face $\gamma_{0} \preceq \gamma$.

\begin{defi}\label{bunchdef}
A {\em bunch\/} in $(E \topto{Q} K, \gamma)$ is
a nonempty collection $\Theta$ of projected faces in $K$ 
with the following property: 
A projected face $\tau_{0} \subset K_{\QQ}$ 
belongs to $\Theta$ if and only if
\begin{equation}\label{bunchcondition}
\emptyset 
\; \ne \; 
\tau_{0}^{\circ} \cap \tau^{\circ} 
\; \ne \; 
\tau^{\circ} 
\qquad \text{holds for all } 
\tau \in \Theta
 \text{ with } \tau \ne \tau_{0}.
\end{equation}
\end{defi}

Let us reformulate this definition in a less formal way.  
We say that two cones $\tau_{1}$ and $\tau_{2}$ {\em overlap}, 
if $\tau_{1}^{\circ} \cap \tau_{2}^{\circ} \neq \emptyset$ holds. 
Now, a nonempty collection $\Theta$ of projected faces 
is a bunch if and only if it has the following properties:
\begin{itemize}
\item any two members of $\Theta$ overlap,
\item there is no pair $\tau_{1}, \tau_{2} \in \Theta$ with $\tau_{1}
\subsetneq \tau_{2}$,
\item if a projected face $\tau_{0}$ overlaps each $\tau \in \Theta$, then
$\tau_{1} \subset \tau_{0}$ for a $\tau_{1} \in \Theta$.
\end{itemize}
 
\begin{exam}\label{weightedprojspace}
Let $E := \ZZ^{n}$, and let $K = \ZZ$. Moreover, fix a sequence
$w_{1}, \ldots, w_{n}$ of positive integers having greatest common
divisor one. This gives an epimorphism
$$ Q \colon E \to K, \qquad e_{i} \mapsto w_{i},$$
where $e_{i}$ denotes the $i$-th canonical base vector. 
Setting $\gamma := \cone(e_{1}, \ldots, e_{n})$, 
we obtain a projected cone 
$(E \topto{Q} K, \gamma)$, and $\Theta := \{Q(\gamma)\}$ is 
a bunch.
\end{exam}

Finally, as the third step, we have to fix the notion 
of a morphism of bunches.
For this, we first have to ``unpack'' the combinatorial information 
contained in a bunch. This is done by constructing a further 
collection of cones:

\begin{defi}
Let  $\Theta$ be a bunch in a projected cone 
$(E \topto{Q} K, \gamma)$. 
The {\em covering collection\/} of $\Theta$ is 
$$ \cov(\Theta) := \{\gamma_{0} \preceq \gamma; \; 
                   \gamma_{0} \text{ minimal with } 
                   Q(\gamma_{0}) \supset \tau \text{ for some } 
                   \tau \in \Theta\}. $$
\end{defi}

As Example~\ref{weightedprojspace} shows, $\cov(\Theta)$ 
will in general comprise much more cones than $\Theta$
itself. 
We can reconstruct the bunch from its covering collection:
$$ 
\Theta 
= 
\{\tau; \; \tau \text{ minimal with } \tau=Q(\gamma_{0})
         \text{ for some } \gamma_{0} \in \cov(\Theta) \}.
$$
In general, for an element
$\gamma_{0} \in \rm \cov(\Theta)$, the image $Q(\gamma_{0})$ need not
be an element of $\Theta$. 
For later purposes, the following observation will be crucial:

\begin{lemma}[Overlapping Property]\label{interior}
Let  $\Theta$ be a bunch in $(E \topto{Q} K, \gamma)$. 
For any two $\gamma_{1}, \gamma_{2} \in
\cov(\Theta)$, we have $Q(\gamma_{1})^{\circ} \cap Q(\gamma_{2})^{\circ}
\ne \emptyset$. 
\end{lemma} 

\proof Let $\sigma_{i} := Q(\gamma_{i})$. By the definition of
$\cov(\Theta)$, there exist cones $\tau_{1}, \tau_{2}
\in \Theta$ with $\tau_{i} \subset \sigma_{i}$. 
Now assume that the relative interiors of the cones
$\sigma_{i}$ are disjoint. Then there is a proper
face $\sigma_{0} \prec \sigma_{1}$ such that $\sigma_{1} \cap
\sigma_{2}$ is contained in $\sigma_{0}$. 

Clearly, $\tau_{1} \cap \tau_{2}$ is contained in
$\sigma_{0}$. Moreover, by the condition~\ref{bunchcondition}, the
intersection $\tau_{1}^{\circ} \cap  
\tau_{2}^{\circ}$ is not empty. In particular, $\tau_{1}^{\circ}$ meets
$\sigma_{0}$. Since $\sigma_{0}$ is a face of $\sigma_{1}$, we
conclude $\tau_{1} \subset \sigma_{0}$. Thus $\gamma_{0} :=
Q^{-1}(\sigma_{0}) \cap \gamma_{1}$ is a proper face of $\gamma_{1}$ 
such that
$Q(\gamma_{0}) = \sigma_{0}$ contains an element of $\Theta$. This
contradicts minimality of $\gamma_{1}$. \endproof 

We come back to the definition of a morphism
of bunches. It is formulated in terms of the 
respective covering collections:

\begin{defi}\label{bunchmorphdef}
Let $\Theta_{i}$ be bunches in projected cones
$(E_{i} \topto{Q_{i}} K_{i}, \gamma_{i})$. 
A {\em morphism\/} from $\Theta_{1}$ to $\Theta_{2}$ is a morphism 
$\Phi \colon E_{1} \to E_{2}$ of the projected cones
such that for every $\alpha_{2} \in \cov(\Theta_{2})$ there is an
$\alpha_{1} \in \cov(\Theta_{1})$ with 
$\Phi(\alpha_{1}) \subset \alpha_{2}$. 
\end{defi}

This concludes the definition of the category of bunches.
The notion of an isomorphism is characterized
as follows:

\begin{prop}\label{bunchiso}
A morphism $\Phi$ of bunches $\Theta_{1}$ and
$\Theta_{2}$ is an isomorphism if and only if 
$\Phi$ is an isomorphism of the ambient projected cones 
and the induced map $\b{\Phi}$ defines a
bijection $\Theta_{1} \to \Theta_{2}$.
\end{prop}

\proof
Let the bunch $\Theta_{i}$ live in the projected cone 
$(E_{i} \topto{Q_{i}} K_{i}, \gamma_{i})$.
Suppose first that $\Phi \colon E_{1} \to E_{2}$ is an
isomorphism of the bunches.
Then there is a morphism of bunches 
$\Psi \colon E_{2} \to E_{1}$ from $\Theta_{2}$ to $\Theta_{1}$
such that $\Phi$ and $\Psi$ are inverse to each other as lattice
homomorphisms.
Note that $\Phi$ and $\Psi$ are as well inverse to each other 
as morphisms of projected cones.

In order to see that  $\b{\Phi} \colon K_{1} \to K_{2}$ defines 
a bijection $\Theta_{1} \to \Theta_{2}$, it suffices to show
that $\Phi$ defines a bijection 
$\cov(\Theta_{1}) \to \cov(\Theta_{2})$.
By bijectivity of $\Phi$ and $\Psi$, we only have to show that
for every $\alpha_{1} \in \cov(\Theta_{1})$ the image
$\Phi(\alpha_{1})$ belongs to $\cov(\Theta_{2})$.
This is done as follows:

Given $\alpha_{1} \in \cov(\Theta_{1})$, we apply 
Definition~\ref{bunchmorphdef} to $\Psi$, and obtain
an $\alpha_{2} \in \cov(\Theta_{2})$ with 
$\Psi(\alpha_{2}) \subset \alpha_{1}$.
Applying $\Phi$ gives $\alpha_{2} \subset \Phi(\alpha_{1})$.
Again by Definition~\ref{bunchmorphdef},
we find an $\t{\alpha}_{1} \in \cov(\Theta_{1})$
with $\Phi(\t{\alpha}_{1}) \subset \alpha_{2}$.
Thus we have $\Phi(\t{\alpha}_{1}) \subset  \Phi(\alpha_{1})$.
By the definition of $\cov(\Theta_{1})$, we obtain
$\t{\alpha}_{1} = \alpha_{1}$.
Consequently, $\Phi(\alpha_{1}) = \alpha_{2}$
belongs to $\cov(\Theta_{2})$.

Now suppose that $\Phi$ is an isomorphism of projected 
cones and that $\b{\Phi}$ defines a bijection 
$\Theta_{1} \to \Theta_{2}$.
Let $\Psi \colon E_{2} \to E_{1}$ denote the inverse
of $\Phi$ as a morphism of projected cones.
The only thing we have to show is that $\Psi$ is 
a morphism from $\Theta_{2}$ to $\Theta_{1}$.
This is done below:

Let $\alpha_{1} \in \cov(\Theta_{1})$. 
Then $\alpha_{2} := \Phi(\alpha_{1})$ is a face of 
$\gamma_{2}$.
We check $\alpha_{2} \in \cov(\Theta_{2})$.
By the definition of $\cov(\Theta_{1})$,
the face $\alpha_{1} \preceq \gamma_{1}$ is minimal 
with $Q_{1}(\alpha_{1}) \supset \tau_{1}$ for some
$\tau_{1} \in \Theta_{1}$.
Thus $\alpha_{2} \preceq \gamma_{2}$ is minimal
with $Q_{2}(\alpha_{2}) \supset \b{\Phi}(\tau_{1})$
for some $\tau_{1} \in \Theta_{1}$.
Since $\b{\Phi}$ induces a bijection 
$\Theta_{1} \to \Theta_{2}$, 
we obtain $\alpha_{2} \in \cov(\Theta_{2})$.
\endproof

The reminder of this section is devoted to the visualization 
of bunches. The idea is that one should be able
to recover many basic properties of a bunch $\Theta$ 
without knowing the ambient projected cone 
$(E \topto{Q} K, \gamma)$.
This will work for the following important class of 
bunches:

\begin{defi}
By a {\em free\/} bunch we mean a bunch $\Theta$ in a 
projected cone $(E \topto{Q} K, \gamma)$, where $\gamma$
is a regular cone in $E$.
\end{defi}

The following construction shows that every free bunch arises
from a certain collections of data in some lattice $K$:

\begin{constr}\label{bunchfromdownstairs}
Let $(w_{1}, \ldots, w_{n})$ be a sequence in 
a lattice $K$ such that the $w_i$ generate $K$.
We speak of the {\em weight vectors} $w_{i}$, and call any 
cone generated by some of the $w_{i}$ a {\em weight cone}.
Let $\Theta$ be a collection of weight cones in $K$ satisfying
Condition~\ref{bunchcondition}, i.e., 
a weight cone $\tau_{0}$ belongs to $\Theta$ if and only if
\begin{equation}
\emptyset 
\; \ne \; 
\tau_{0}^{\circ} \cap \tau^{\circ} 
\; \ne \; 
\tau^{\circ} 
\qquad \text{holds for all } 
\tau \in \Theta
 \text{ with } \tau \ne \tau_{0}.
\end{equation}
Then there is an associated projected cone $(E \topto{Q} K, \gamma)$
with the lattice $E := \ZZ^{n}$, the cone 
$\gamma := \cone(e_{1}, \ldots, e_{n})$ spanned by the canonical 
base vectors, and the map $Q \colon E \to K$ sending $e_{i}$ to $w_{i}$.
By construction, the collection $\Theta$ is a bunch in 
the projected cone $(E \topto{Q} K, \gamma)$.
\end{constr}

This construction allows us to visualize bunches.
We regard a given sequence of weight vectors as 
a set $\{w_{1}, \dots, w_{n}\}$, where each $w_{i}$ has a
multiplicity $\mu_{i}$ counting the number of its
repetitions. 
Then we may put these data as well as the cones of 
a given bunch in a picture.
For example, in the setting of~\ref{weightedprojspace}, 
the bunch $\Theta$ defined by the
sequence $(1,3,5,5)$ arises from 
the picture
\begin{center}
  \begin{picture}(0,0)%
\includegraphics{weightedprojspce.pstex}%
\end{picture}%
\setlength{\unitlength}{1243sp}%
\begingroup\makeatletter\ifx\SetFigFont\undefined%
\gdef\SetFigFont#1#2#3#4#5{%
  \reset@font\fontsize{#1}{#2pt}%
  \fontfamily{#3}\fontseries{#4}\fontshape{#5}%
  \selectfont}%
\fi\endgroup%
\begin{picture}(6356,767)(439,-138)
\put(4726,389){\makebox(0,0)[lb]{\smash{\SetFigFont{6}{7.2}{\familydefault}{\mddefault}{\updefault}
$(1)$
}}}
\put(5626,389){\makebox(0,0)[lb]{\smash{\SetFigFont{6}{7.2}{\familydefault}{\mddefault}{\updefault}
$(2)$
}}}
\put(3826,389){\makebox(0,0)[lb]{\smash{\SetFigFont{6}{7.2}{\familydefault}{\mddefault}{\updefault}
$(1)$
}}}
\end{picture}

\end{center}

As one might expect, this picture describes the three-dimensional
weighted projective space $\PP_{1,3,5,5}$.
Moreover, as we shall see in Proposition~\ref{Kleinschmidt},
the smooth complete toric varieties with Picard group 
$\ZZ^{2}$ arise from sequences in $\ZZ^{2}$ and 
a collection $\Theta = \{\tau\}$ according to the following
picture:

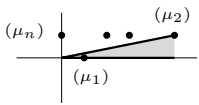
\begin{center}
  \begin{picture}(0,0)%
\includegraphics{smoothclass.pstex}%
\end{picture}%
\setlength{\unitlength}{1243sp}%
\begingroup\makeatletter\ifx\SetFigFont\undefined%
\gdef\SetFigFont#1#2#3#4#5{%
  \reset@font\fontsize{#1}{#2pt}%
  \fontfamily{#3}\fontseries{#4}\fontshape{#5}%
  \selectfont}%
\fi\endgroup%
\begin{picture}(3847,1937)(1126,-2333)
\put(2476,-1861){\makebox(0,0)[lb]{\smash{\SetFigFont{6}{7.2}{\familydefault}{\mddefault}{\updefault}
$(\mu_{1})$
}}}
\put(1126,-961){\makebox(0,0)[lb]{\smash{\SetFigFont{6}{7.2}{\familydefault}{\mddefault}{\updefault}
$(\mu_{n})$
}}}
\put(4101,-636){\makebox(0,0)[lb]{\smash{\SetFigFont{6}{7.2}{\familydefault}{\mddefault}{\updefault}
$(\mu_{2})$
}}}
\end{picture}

\end{center}

In order to compare two bunches arising from 
Construction~\ref{bunchfromdownstairs}, 
there is no need to determine the covering collection. 
Namely, using Proposition~\ref{bunchiso}, we obtain:

\begin{rem}\label{downstairsisochar}
Two sets of data $(K;w_{1}, \ldots, w_{n}; \Theta)$
and  $(K';w_{1}', \ldots, w_{n}'; \Theta')$ 
as in~\ref{bunchfromdownstairs} have isomorphic
associated free bunches if and only if there is a lattice 
isomorphism $\b{\Phi} \colon K \to K'$ such that
\begin{enumerate}
\item $(w_{1}', \ldots, w_{n}')$ and
  $(\b{\Phi}(w_{1}), \ldots, \b{\Phi}(w_{n}))$
  differ only by enumeration,
\item the collections $\Theta'$ and $\{\b{\Phi}(\tau); \; \tau \in
  \Theta\}$ coincide. 
\end{enumerate}
\end{rem}

\section{The basic duality lemmas}
\label{dualityinweightedlattices}

In this section, we provide basic duality 
statements for translating from the language
of bunches into the language of fans.
First we need a concept of a dual of a given projected cone
$(E \topto{Q} K, \gamma)$.
For this, note that $(E \topto{Q} K, \gamma)$ determines 
two exact sequences of lattice homomorphisms
$$ \xymatrix{
    0 \ar[r] &
    {M} \ar[r] &
    {E} \ar[r]^{Q} & 
    K \ar[r] &  
    0,  }$$
$$ \xymatrix{
    0 \ar[r] &
    {L} \ar[r] &
    {F} \ar[r]^{P} & 
    {N} \ar[r] &  
    0 , }$$
where $M$ is the kernel of $Q \colon E \to K$, and the 
second sequence arises from the first one by 
applying $\Hom(?,\ZZ)$; note that
$P$ is not the dual homomorphism of $Q$.
Let $\delta := \gamma^{\vee}$ denote the dual cone. Then $\delta$
is again strictly convex, simplicial and of full dimension.

\begin{defi}
We call $(F \topto{P} N, \delta)$ the 
{\em dual projected cone\/} of 
$(E \topto{Q} K, \gamma)$.
\end{defi}

In the sequel, fix a projected cone $(E \topto{Q} K, \gamma)$,
and denote the associated dual projected cone by 
$(F \topto{P} N, \delta)$. Recall that we have 
the face correspondence, see for 
example~\cite[Appendix~A]{Od}:

\begin{rem}
The sets of faces of the cones $\gamma \subset E_{\QQ}$ and 
$\delta \subset F_{\QQ}$ are in order reversing correspondence via
$$ \faces(\gamma)  \to \faces(\delta),
\qquad
\gamma_{0}  \mapsto  \gamma_{0}^{*} := \gamma_{0}^{\perp} \cap \delta.
$$
\end{rem}

Our task is to understand the relations
between the projected faces $Q(\gamma_{0})$,
where $\gamma_{0} \preceq \gamma$,
and the images $P(\gamma_{0}^{*})$
of the corresponding faces. 
The following observation is central:

\begin{lemma}[Invariant Separation Lemma]%
\label{invariantseparation}
Let $\gamma_{1}, \gamma_{2} \preceq \gamma$, let $\delta_{i} :=
\gamma_{i}^{*}$, and let $L := \ker(P)$. Then the following 
statements are equivalent:   
\begin{enumerate}
\item There is an $L$-invariant separating linear form for
  $\delta_{1}$ and $\delta_{2}$.
\item For the relative interiors $Q(\gamma_{i})^{\circ}$ we have
  $Q(\gamma_{1})^{\circ} \cap Q(\gamma_{2})^{\circ} \ne \emptyset$.
\end{enumerate}
\end{lemma}

Here, by an {\em $L$-invariant\/} linear form we mean 
an element $u \in E = \Hom(F,\ZZ)$ with $L \subset u^{\perp}$. 
Moreover, recall from Section~\ref{section1} that a separating 
linear form for the cones $\delta_{1}$ and $\delta_{2}$ is an 
element $u \in E$ with
$$
u\vert_{\delta_{1}} \ge 0,
\qquad
u\vert_{\delta_{2}} \le 0,
\qquad
u^{\perp} \cap \delta_{1} 
= u^{\perp} \cap \delta_{2} 
= \delta_{1} \cap \delta_{2}.$$ 

\proof[Proof of Lemma~\ref{invariantseparation}]
As before, let $M := \ker(Q)$. Then the $L$-invariant
linear forms on $F$ are precisely the elements of $M$.
Thus, since $\delta_1 \cap \delta_2$ is a face of both
$\delta_i$, condition~(i) is equivalent to
\begin{eqnarray}
\label{invsepchar}
M 
\cap 
\left(\delta_{1}^{\vee} \cap (\delta_{1} \cap \delta_{2})^{\perp}
\right)^{\circ} 
\cap 
- \left(\delta_{2}^{\vee} \cap (\delta_{1} \cap \delta_{2})^{\perp}
\right)^{\circ} 
& \ne & 
\emptyset.
\end{eqnarray}

Note that $\delta_{i}^{\vee}$ equals $\lin(\gamma_{i}) +
\gamma$. Moreover, $(\delta_{1} \cap \delta_{2})^{\perp}$ equals
$\lin(\gamma_{1}) + \lin(\gamma_{2})$, because the cone $\delta$
is simplicial. Hence we obtain
\begin{eqnarray*}
\left(\delta_{1}^{\vee} \cap (\delta_{1} \cap \delta_{2})^{\perp}
\right)^{\circ} 
& = & 
\left((\lin(\gamma_{1}) + \gamma) 
    \cap (\lin(\gamma_{1}) +  \lin(\gamma_{2})) \right)^{\circ} \\
& = & 
\lin(\gamma_{1}) + \gamma_{2}^{\circ}. 
\end{eqnarray*}

For the second equality we used simpliciality of $\gamma$. 
By analogous arguments, the expression
$-(\delta_{2}^{\vee} \cap (\delta_{1} \cap
\delta_{2})^{\perp})^{\circ}$
simplifies to $\lin(\gamma_{2}) - \gamma_{1}^{\circ}$.
Thus~\ref{invsepchar} holds if and only if 
\begin{eqnarray}
\label{invsepchar2}
M 
\cap 
\left(\lin(\gamma_{1}) + \gamma_{2}^{\circ}\right) 
\cap 
\left(\lin(\gamma_{2}) - \gamma_{1}^{\circ}\right) 
& \ne & 
\emptyset.
\end{eqnarray}

We claim that the left hand side simplifies to 
$M \cap (\gamma_{2}^{\circ} - \gamma_{1}^{\circ})$. 
Indeed, any $u \in M$ belonging to the left hand side
of~\ref{invsepchar2} has a unique representation 
in terms of the primitive generators of 
$e_{1}, \ldots, e_{n}$ of $\gamma$:
$$
u 
\; = \; 
\sum_{e_{i} \in \gamma_{1} \setminus \gamma_{2}} a_{i} e_{i}
+
\sum_{e_{j} \in \gamma_{1} \cap \gamma_{2}} b_{j} e_{j}
+
\sum_{e_{k} \in \gamma_{2} \setminus \gamma_{1}} c_{l} e_{l},
\qquad
\text{where } a_{i} < 0 \text{ and } c_{l} > 0.
$$

Now, dividing the middle term into two sums, one with
only positive coefficients and the other with only 
negative ones, gives 
$u \in \gamma_{2}^{\circ} - \gamma_{1}^{\circ}$.
The reverse inclusion is obvious.

Consequently~\ref{invsepchar2} is equivalent to 
$M \cap (\gamma_{2}^{\circ} - \gamma_{1}^{\circ}) \ne \emptyset$.
This in turn is obviously equivalent to condition (ii).
\endproof
 
Let us mention here that simpliciality of the cones $\gamma$ and
$\delta$ is essential for the Invariant Separation Lemma:

\begin{exam}
Consider the ``projected cone'' 
$(E \topto{Q} K, \gamma)$,
where the lattices are $E := \ZZ^{3}$ and $K := \ZZ$,
the map $Q$ is the projection onto the third coordinate,
and the cone $\gamma$ is given in terms of canonical 
base vectors by
$$ 
\gamma = 
\cone 
(e_{1}+e_{3},e_{2}+2e_{3},e_{1}-2e_{3},e_{2}-e_{3}). 
$$
Denote the ``dual projected cone'' by  
$(F \topto{P} N, \delta)$.
Then $L := \ker(P)$ is the sublattice generated by the dual base 
vector $e_{3}^{*}$, and the cone $\delta$ is given by
$$
\delta 
= 
\cone 
(e_{1}^{*},e_{2}^{*},
 e_{1}^{*}+2e_{2}^{*}-e_{3}^{*},
 2e_{1}^{*}+e_{2}^{*}+e_{3}^{*}). 
$$
The faces $\gamma_{1} := \cone(e_{2}-e_{3})$ and 
$\gamma_{2} := \cone(e_{1}+e_{3})$ do not satisfy
Condition~\ref{invariantseparation}~(ii).
Nevertheless, the corresponding faces
$$
\gamma_{1}^{*} = \cone(e_{1}^{*},2e_{1}^{*}+e_{2}^{*}+e_{3}^{*}),
\qquad 
\gamma_{2}^{*} = \cone(e_{2}^{*},e_{1}^{*}+2e_{2}^{*}-e_{3}^{*}),
$$
admit $L$-invariant separating linear forms. For example
we can take the linear form $e_{1}-e_{2} \in E$.
\end{exam}

Next we compare injectivity of $Q$ with surjectivity of $P$ along
corresponding faces (of course, the roles of $Q$, $\delta_{0}$ etc.
and $P$, $\gamma_{0}$, etc. can be interchanged in the statement): 

\begin{lemma}\label{surjectiveandinjective1}
For a face $\gamma_{0} \preceq \gamma$ and
$\delta_{0} := \gamma_{0}^{*}$,
the following statements are equivalent:
\begin{enumerate}
\item $P$ maps $\lin(\delta_{0})$ onto $N_{\QQ}$,
\item $Q$ is injective on $\lin(\gamma_{0})$.
\end{enumerate}
\end{lemma}

\proof Let $M := \ker(Q)$ and $L := \ker(P)$. Using the fact that
$\lin(\delta_{0})$ and $\lin(\gamma_{0})$ are the orthogonal spaces of
each other, we obtain the assertion by dualizing:
\begin{eqnarray*}
\lin(\delta_{0}) + L_{\QQ} = F_{\QQ} 
& \iff &
\lin(\gamma_{0}) \cap M_{\QQ} = \{0\}. \qquad \qed
\end{eqnarray*}

If we take the lattice structure into consideration,
then the situation becomes slightly more involved. 
The essential observation is: 

\begin{lemma}\label{surjectiveandinjective2}
For a face $\gamma_{0} \preceq \gamma$ and
$\delta_{0} := \gamma_{0}^{*}$,
the following statements are equivalent:
\begin{enumerate}
\item $P$ maps $\lin(\delta_{0}) \cap F$ onto $N$,
\item $Q$ maps $\lin(\gamma_{0}) \cap E$ isomorphically onto a
  primitive sublattice of $K$.
\end{enumerate}
\end{lemma}

\proof Set $L := \ker(P)$. Assume that (i) holds. Then the snake lemma
provides an exact sequence
$$ \xymatrix{ 
    0 \ar[r] & 
    {L \cap \lin(\delta_{0})} \ar[r] & 
    {L} \ar[r] &
    {F/(\lin(\delta_{0}) \cap F)} \ar[r] & 
    0. } $$ 

The dual lattice of $F/(\lin(\delta_{0}) \cap F)$ is canonically
isomorphic to $E \cap \lin (\gamma_{0})$. Hence, applying
$\Hom(?,\ZZ)$ gives an exact sequence 
$$ \xymatrix{
  0 \ar[r] & 
  E \cap \lin(\gamma_{0}) \ar[r]^(.625){Q} & 
  K \ar[r] & 
  \Hom(L \cap \lin(\delta_{0}),\ZZ) \ar[r] & 
  0.} $$

This implies condition~(ii). The reverse direction can be settled by
similar arguments. \endproof

Finally, we consider morphisms $\Phi \colon E_{1} \to E_{2}$ of
projected cones $(E_{i} \topto{Q_{i}} K_{i}, \gamma_{i})$.  
These data define a commutative diagram of lattices with exact rows:
$$ 
\xymatrix{
    0 \ar[r] &
    {M_{1}} \ar[r]  \ar[d]&
    {E_{1}} \ar[r]^{Q_{1}} \ar[d]_{{\Phi}}& 
    {K_{1}} \ar[r] \ar[d]^{{\bar{\Phi}}}&  
    0 \\
    0 \ar[r] &
    {M_{2}} \ar[r] &
    {E_{2}} \ar[r]_{Q_{2}}  & 
    {K_{2}} \ar[r]  &  
    0  }
$$

Applying $\Hom(?,\ZZ)$ to this diagram, 
we obtain the following commutative diagram,
again with exact rows:
$$ 
\xymatrix{
    0 \ar[r] &
    {L_{1}} \ar[r]  &
    {F_{1}} \ar[r]^{P_{1}} & 
    {N_{1}} \ar[r] &  
    0 \\
    0 \ar[r] &
    {L_{2}} \ar[r] \ar[u] &
    {F_{2}} \ar[r]_{P_{2}}  \ar[u]^{{\Psi}}& 
    {N_{2}} \ar[r]  \ar[u]_{{\bar{\Psi}}}&  
    0  }
$$

\begin{rem}
The dual map ${\Psi} \colon F_{2} \to F_{1}$ is
a morphism of the dual projected cones 
$(F_{i} \topto{P_{i}} N_{i}, \delta_{i})$
of 
$(E_{i} \topto{Q_{i}} K_{i}, \gamma_{i})$.
\end{rem}

Now, consider faces $\alpha_{i} \preceq \gamma_{i}$, and let $\beta_{i} :=
\alpha_{i}^{*}$ denote the corresponding faces of the cones $\delta_{i}$. 

\begin{lemma}\label{shiftingmorphisms}
We have 
$\Phi(\alpha_{1}) \subset \alpha_{2}$ if and only if 
$\Psi (\beta_{2}) \subset \beta_{1}$ holds.
\end{lemma}

\proof 
We only have to verify one implication. The other then is a simple
consequence of $\alpha_{i} = \beta_{i}^{*}$. So, suppose
$\Phi (\alpha_{1}) \subset \alpha_{2}$. Then we obtain
$$ 
\Psi (\beta_{2}) 
\; = \; 
{\Psi} (\alpha_{2}^{\perp} \cap \delta_{2})
\;  \subset \; 
\Psi (\alpha_{2}^{\perp}) \cap {\Psi} (\delta_{2}) 
\; \subset \; 
\alpha_{1}^{\perp} \cap \delta_{1} 
\; = \; \beta_{1}.
\qquad \qed
$$

\section{Bunches and fans}
\label{bunchesandadmissiblefans}

In this section, we compare bunches with fans.
We shall show that the category of bunches is equivalent to 
the category of ``maximal projectable fans'',
see Theorem~\ref{dualitythm}. 
The latter category is defined as follows: 

\begin{defi}\label{admissiblefandefi}
\begin{enumerate}
\item Let $(F \topto{P} N, \delta)$ be a projected cone, and 
let $L := \ker(P)$.
\begin{enumerate}
\item A {\em projectable fan\/} in $(F \topto{P} N, \delta)$ 
  is a fan  $\Sigma$ consisting of faces of $\delta$ such that
  any two maximal cones of $\Sigma$ can be separated by an
  $L$-invariant linear form.
\item We call a projectable fan $\Sigma$ in $(F \topto{P} N, \delta)$ {\em
maximal\/} if any $\delta_{0} \prec \delta$, which can be separated by
$L$-invariant linear forms from the maximal cones of $\Sigma$, belongs to
$\Sigma$.
\end{enumerate}
\item A {\em morphism\/} of projectable fans $\Sigma_{i}$ 
(maximal or not) 
in projected cones $(F_{i} \topto{P_{i}} N_{i}, \delta_{i})$ is a
morphism $\Psi \colon F_{1} \to F_{2}$ of
projected cones which is in addition a map of the fans
$\Sigma_{i}$. 
\end{enumerate}
\end{defi}

Note that a projectable fan is the collection of all faces
of the cones belonging to a ``locally coherent costring''
in the sense of~\cite[Def.~5.1]{Hu}, but the converse does 
not hold in general.
We shall demonstrate later by means of an example
the importance of the maximality 
condition~(b), see~\ref{underlinemaximal}.

We define now a functor $\mathfrak{F}$ from bunches to maximal projectable
fans. Let $\Theta$ be a bunch in the projected cone 
$(E \topto{Q} K, \gamma)$. 
Consider the associated dual projected cone
$(F \topto{P} N, \delta)$ and
the following subfan of the fan of faces of $\delta$:
$$ \Sigma := \{\sigma \preceq \delta; \; \sigma \preceq
\gamma_{0}^{*} \text{ for some } \gamma_{0} \in \cov(\Theta)\}. $$

\begin{lemma}
$\Sigma$ is a maximal projectable fan in 
$(F \topto{P} N, \delta)$.
\end{lemma}

\proof Let $L := \ker(F)$. By the Overlapping Property~\ref{interior}
of the Covering Collection $\cov(\Theta)$ and the Invariant Separation
Lemma~\ref{invariantseparation}, any two maximal cones of $\Sigma$
can be separated by $L$-invariant linear forms. So we only have to
verify the maximality condition~\ref{admissiblefandefi}~(i)~(b) 
for $\Sigma$. 

Suppose that the face $\sigma \preceq \delta$ can be separated by
$L$-invariant linear forms from the maximal cones of $\Sigma$ but
does not belong to $\Sigma$. 
The projected face $\tau_{0} := Q(\sigma^{*})$ does not belong to
$\Theta$, because otherwise any minimal face 
$\gamma_{0} \preceq \sigma^{*}$ projecting onto $\tau_{0}$
would belong to $\cov(\Theta)$, which contradicts the 
choice of $\sigma$.

The Invariant Separation Lemma~\ref{invariantseparation} yields 
$\tau_{0}^{\circ} \cap \tau^{\circ} \ne \emptyset$ for every $\tau \in
\Theta$. Since $\tau_{0}$ is not an element of $\Theta$, it has to
contain some element of $\Theta$. But then some face of 
$\sigma^{*}$ belongs to the Covering Collection $\cov(\Theta)$.
Again this contradicts the choice of the face $\sigma \preceq
\delta$. \endproof

The assignment $\Theta \mapsto \Sigma$
extends canonically to morphisms. Namely,
let $\Theta_{i}$ be bunches in projected cones 
$(E_{i} \topto{Q_{i}} K_{i}, \gamma_{i})$. 
Let $\Sigma_{i}$ denote the associated maximal projectable fans in 
the respective dual projected cones  
$(F_{i} \topto{P_{i}} N_{i}, \delta_{i})$. 

Lemma~\ref{shiftingmorphisms} tells us that for every morphism
$\Phi \colon E_{1} \to E_{2}$ of the bunches $\Theta_{1}$ and
$\Theta_{2}$, the dual map $\Psi \colon F_{2} \to F_{1}$ is a
morphism of the maximal projectable fans $\Sigma_{2}$ and $\Sigma_{1}$.
Thus, we obtain:

\begin{prop}
The assignments $\Theta \mapsto \Sigma$ and $\Phi \mapsto
\Psi$ define a contravariant functor $\mathfrak{F}$ from the
category of bunches to the category of maximal projectable fans. 
\endproof
\end{prop}

Now we go the other way round. Consider a maximal projectable fan 
$\Sigma$ in a projected cone 
$(F \topto{P} N, \delta)$. 
Let $(E \topto{Q} K, \gamma)$
denote the associated dual projected cone.
Define $\Theta$ to be the set of the minimal cones among all projected
faces in $K$ arising from $\Sigma$:
$$ \Theta := \{\tau_{0}; \; \tau_{0} \text{ minimal with } \tau_{0} =
Q(\delta_{0}^{*}) \text{ for some } \delta_{0} \in \Sigma^{\max}\}. $$

\begin{lemma}
$\Theta$ is a bunch in $(E \topto{Q} K, \gamma)$.
\end{lemma}

\proof We verify Property~\ref{bunchcondition} for a given
$\tau_{0} \in \Theta$. According to the Invariant Separation
Lemma~\ref{invariantseparation}, we have 
$\tau_{0}^{\circ} \cap \tau^{\circ} \ne \emptyset$ 
for any further $\tau \in \Theta$. Moreover, since
$\Theta$ consists of minimal cones, $\tau^{\circ}$ is not contained
in $\tau_{0}^{\circ}$ provided that $\tau \in \Theta$ is 
different from $\tau_{0}$.

Conversely, let the projected face $\tau_{0}$
satisfy~\ref{bunchcondition}. Choose $\gamma_{0} \preceq \gamma$
with $Q(\gamma_{0}) = \tau_{0}$. The Invariant Separation
Lemma~\ref{invariantseparation} tells us that $\delta_{0} := \gamma_{0}^{*}$
belongs to $\Sigma$.  Let $\delta_{1} \in \Sigma$ be a maximal cone with
$\delta_{0} \preceq \delta_{1}$, and consider the image $\tau_{1} :=
Q(\delta_{1}^{*})$.  Then we have $\tau_{1} \subset \tau_{0}$, because
$\delta_{1}^{*} \preceq \delta_{0}^{*}$ holds.

By the definition of the collection $\Theta$, there is a cone
$\tau_{2} \in \Theta$ with $\tau_{2} \subset \tau_{1}$.
In particular, we have $\tau_{2} \subset \tau_{0}$.
Applying once more  the Invariant Separation 
Lemma~\ref{invariantseparation},
gives even $\tau_{2}^{\circ} \subset \tau_{0}^{\circ}$.
Thus Property~\ref{bunchcondition} yields $\tau_{0}
= \tau_{2}$. This shows $\tau_{0} \in \Theta$. \endproof

According to Lemma~\ref{shiftingmorphisms},
associating to a map $\Psi$ of maximal projectable fans its dual map
$\Phi$ makes this construction functorial. Thus we have:

\begin{prop}\label{projectablefan2bunch}
The assignments $\Sigma \mapsto \Theta$ and $\Psi \mapsto
\Phi$ define a contravariant functor $\mathfrak{B}$ from the
category of maximal projectable fans to the category of bunches. \endproof 
\end{prop}

Summing up, we arrive at the main result of this section, namely the
following duality statement:

\begin{thm}\label{dualitythm}
The functors $\mathfrak{F}$ and $\mathfrak{B}$ are inverse to each
other. In particular, the categories of bunches and maximal projectable fans
are dual to each other. \endproof
\end{thm}

Let us emphasize here the role of the maximality 
condition~\ref{admissiblefandefi}~(i)~(b) 
in this result. The following example shows that there
is no hope for a similar statement on (non maximal)
projectable fans:

\begin{exam}\label{underlinemaximal}
Consider the projected cone $(F \topto{P} N, \delta)$, where
the lattices are $F := \ZZ^{3}$ and $N := \ZZ^{2}$, the cone 
$\delta$ is the positive orthant in $F_{\QQ}$, and the projection
map is given by
$$ 
P \colon F \to N,
\qquad
(v_{1},v_{2},v_{3}) \mapsto (v_{1}-v_{3},v_{2}-v_{3}).
$$

Then the fan $\Sigma$ in $F$ having 
$\delta_{1} := \cone(e_{1},e_{2})$
and 
$\delta_{2} := \cone(e_{1},e_{3})$ 
as its maximal cones is projectable. 
But $\Sigma$ is not maximal, because we may
enlarge it to a (maximal) projectable fan
$\Sigma'$ by adding the cone
$\delta_{3} := \cone(e_{2},e_{3})$.

Let $(E \topto{Q} K, \gamma)$ be the dual projected cone of 
$(F \topto{P} N, \delta)$. 
Then we have $E \cong \ZZ^{3}$ and $K \cong \ZZ$.
Moreover the projection $Q$ sends each dual base vector
$e^{*}_{i}$ to $1 \in \ZZ$. In particular, we obtain
$$ 
Q(\delta_{i}^{*})
\; = \; 
Q(\delta_{i}^{\perp} \cap \gamma) 
\; = \; 
\QQ_{\ge 0}, 
\quad i = 1,2,3.
$$ 

Thus, $\Sigma$ and $\Sigma'$ determine the same collection
$\Theta = \{\QQ_{\ge 0}\}$ of projected faces in $K$.
In other words, there is no way to reconstruct $\Sigma$ via 
face duality from a collection of projected faces in $K$.
\end{exam}

In the rest of this section, we associate to any
maximal projectable fan its ``quotient fan''. 
So, let $\Sigma$ be a maximal projectable fan in a weighted
lattice $(F \topto{P} N, \delta)$. Then the images
$P(\sigma)$, where $\sigma$ runs through the maximal cones
of $\Sigma$ are the maximal cones of a quasifan $\Sigma'$ in
$N$.
 
We reduce $\Sigma'$ to a fan as follows: Let $L' \subset N$ be the
primitive sublattice generating the minimal cone of $\Sigma'$, let
$N' := N/L'$, and let $P' \colon N \to N'$ denote the projection.

\begin{defi}
The {\em quotient fan\/} of $\Sigma$ is 
$\Delta := \{P'(\sigma'); \; \sigma' \in \Sigma'\}$.
\end{defi}

Note that $R := P' \circ P \colon F \to N'$ is a map of the
fans $\Sigma$ and $\Delta$. In fact, this is a special case of a
more general construction, see~\cite[Theorem~2.3]{AcHa}. In our
setting, it is easy to see that everything is compatible with
morphisms. Thus we obtain:  

\begin{prop}\label{quotfan}
The assignment $\Sigma \mapsto \Delta$ defines a covariant
functor $\mathfrak{Q}$ from the category of maximal projectable 
fans to the category of fans. \endproof
\end{prop}

The following simple example shows that dividing by $L'$ in
the construction of quotient fan is indeed necessary:

\begin{exam}
Consider the lattices $E := \ZZ^{2}$ and $K := \ZZ$, the map 
$Q \colon E \to K$, $(u_{1},u_{2}) \mapsto u_{1} + u_{2}$, and the
positive orthant $\gamma \subset \QQ^{2}$. Let $\Theta$ be the bunch
consisting just of the trivial cone $\{0\}$. Then the quasifan
$\Sigma'$ determined by $\Theta$ consists of the single cone
$\sigma := \QQ$.
\end{exam} 

We note an observation on the composition 
$\mathfrak{Q} \circ \mathfrak{F}$. Consider a bunch
$\Theta$ in $(E \topto{Q} K, \gamma)$ and its
associated maximal projectable fan $\Sigma$ in 
$(F \topto{P} N, \delta)$. Let
$\Delta$ be the quotient fan of $\Sigma$, and, as before, let
$R \colon F \to N'$ be the projection.

\begin{prop}\label{conesofquotfan}
There is a canonical order reversing bijection
$$ 
\{\gamma_{0} \preceq \gamma; \; \tau^{\circ} \subset Q(\gamma_{0})^{\circ}
                                \text{ for some } \tau \in \Theta\} 
\to \Delta,
\qquad 
\gamma_{0} \mapsto R (\gamma_{0}^{*}).
$$
\end{prop}

\proof The inverse map is given by $\sigma \mapsto (R^{-1}(\sigma)
\cap |\Sigma|)^{*}$. \endproof

\section{Combinatorics of quotients}
\label{Tmaximalsubsets}

Here we present the first application of the language of bunches.
We consider the action of a subtorus on a $\QQ$-factorial nondegenerate
affine toric variety and give a combinatorial description of the
maximal open subsets admitting a good quotient by this action. This
complements results of~\cite{BBSw} for torus actions 
on $X = \CC^{n}$. 

Let us first recall the basic concepts concerning good quotients. 
Let the reductive group $G$ act on a variety $X$ by means of a 
morphism $G \times X \to X$. A {\em good quotient\/} for this
action is a $G$-invariant affine morphism $p \colon X \to Y$ such
that the canonical map $\mathcal{O}_{Y} \to
p_{*}(\mathcal{O}_{X})^{G}$ is an isomorphism. 
If it exists, then the good quotient space is usually 
denoted by $X \quot G$. 

In general, a $G$-variety $X$ need not admit a good
quotient $X \to X \quot G$, 
but there frequently exist many invariant open subsets
$U \subset X$ with good quotient $U \to U \quot G$. 
It is one of the central tasks of Geometric Invariant Theory to
describe all these open subsets, see~\cite[Section~7.2]{BB}.
In the course of this problem, one reasonably looks for maximal
$U \subset X$ in the following sense, see~\cite[Section~7.2]{BB}
and~\cite{BBSw}:

\begin{defi}
An open subset $U \subset X$ is called {\em $G$-maximal}, if there is
a good quotient $p \colon U \to U \quot G$ and there is no open $U'
\subset X$ admitting a good quotient $p' \colon U' \to U' \quot G$
such that $U$ is a proper $p'$-saturated subset of $U'$.
\end{defi}
 
In the setting of subtorus actions, the maximal open subsets with 
good quotient can be characterized in terms of fans. 
This relies on the following observation due to \'{S}wi\c{e}cicka,
see~\cite[Proposition~2.5]{Sw1}:

\begin{prop}\label{swiecicka}
Let $X$ be a toric variety, and let $T \subset T_{X}$ be a subtorus of
the big torus. If $U \subset X$ is a $T$-maximal subset, then $U$ is 
invariant under $T_{X}$.
\end{prop}

Now, let $X$ be the toric variety arising from a fan
$\Delta$ in a lattice $N$, let $T \subset T_{X}$ be the subtorus
corresponding to a primitive sublattice $L \subset N$. 
By the above proposition, the $T$-maximal subsets $U \subset X$
correspond to certain subfans of $\Delta$.
The characterization of these fans is standard, see 
e.g.~\cite[Proposition~1.3]{Ha1}:

\begin{prop}\label{goodquotchar}
Let $U \subset X$ be the open $T_{X}$-invariant subset defined by a
subfan $\Sigma$ of $\Delta$.
\begin{enumerate}
\item There is a good quotient $U \to U \quot T$ if and only if
  any two maximal cones of $\Sigma$ can be separated by an
  $L$-invariant linear form.
\item $U$ is $T$-maximal if and only if (i) holds and every $\sigma
  \in \Delta$ that can be separated by $L$-invariant linear forms from
  the maximal cones of $\Sigma$ belongs to $\Sigma$. \endproof
\end{enumerate}
\end{prop}

Though this is a complete combinatorial description of all 
$T$-maximal subsets, it has two drawbacks in practice: 
On the one hand, the ambient space of the combinatorial
data might be of quite big dimension, and, 
on the other hand, for the explicit checking of the 
conditions there may be large numbers of cones to go through.
The language of bunches makes the situation more clear.

Let $X$ be an affine nondegenerate $\QQ$-factorial toric variety
arising from a cone $\delta$ in a lattice $F$, and let
$T \subset T_{X}$ be the subtorus corresponding to a sublattice
$L \subset F$.
Setting $N := F/L$, we obtain a projected cone
$(F \topto{P} N, \delta)$. 
Moreover, we have the dual projected cone 
$(E \topto{Q} K, \gamma)$,
where $K$ is canonically isomorphic to the lattice of characters 
of the small torus $T \subset T_{X}$.

In order to describe the $T$-maximal subsets of $X$,
we use the functor 
$\mathfrak{F}$ associating to a bunch $\Theta$ in 
$(E \topto{Q} K, \gamma)$ a maximal projectable fan 
$\mathfrak{F}(\Theta)$ in 
$(F \topto{P} N, \delta)$.
The resulting statement generalizes and complements the 
results of~\cite{BBSw}: 

\begin{thm}\label{bunches2Tmax}
The assignment $\Theta \mapsto X_{\mathfrak{F}(\Theta)}$ defines a
one-to-one correspondence between the bunches in 
$(E \topto{Q} K, \gamma)$ and the $T$-maximal open subsets of $X$.
\end{thm}

\proof 
By the definition of a maximal projectable fan and 
Proposition~\ref{goodquotchar}, the
toric open subvariety $X_{\mathfrak{F}(\Theta)}$ is indeed $T$-maximal. Hence
the assignment is well defined. Moreover, it is of course injective.
Surjectivity follows from Proposition~\ref{swiecicka}.
\endproof

\begin{rem}
In the setting of Theorem~\ref{bunches2Tmax}, the good quotient
of the $T$-action on  $X_{\mathfrak{F}(\Theta)}$ is the toric morphism
$X_{\mathfrak{F}(\Theta)} \to  X_{\mathfrak{Q}(\mathfrak{F}(\Theta))}$
arising from the projection 
$\mathfrak{F}(\Theta) \to \mathfrak{Q}(\mathfrak{F}(\Theta))$
onto the quotient fan.
\end{rem}

Some of the good quotients are of special interest:
A {\em geometric quotient\/} for an action of a reductive 
group $G$ on a variety $X$ is a good quotient
that separates orbits. 
Geometric quotients are denoted by $X \to X/G$.

Again, for subtorus actions on toric varieties,
there is a description in terms of fans. Let $X$ be the toric variety
arising from a fan $\Delta$ in a lattice $N$, let $T \subset T_{X}$ be
the subtorus corresponding to a primitive sublattice $L \subset N$.
Existence of a geometric quotient is characterized as follows,
see~\cite[Theorem~5.1]{Hm}: 

\begin{prop}
The action of $T$ on $X$ admits a geometric quotient if and only if
the projection $P \colon N \to N/L$ is injective on the support
$\Sigma$. 
\end{prop}

Let us translate this into the language of bunches.
As before, consider a projected cone 
$(F \topto{P} N, \delta)$ and its associated 
dual projected cone 
$(E \topto{Q} K, \gamma)$.
 
\begin{defi}
A bunch $\Theta$ in $(E \topto{Q} K, \gamma)$ is
called {\em geometric\/} if $\dim(\tau) = \dim(K)$ holds for
every $\tau \in \Theta$. 
\end{defi}

We consider the affine toric variety $X := X_{\delta}$ 
and the subtorus $T \subset T_{X}$ corresponding to the sublattice $L
\subset N$. The above notion yields what we are looking for:

\begin{prop}\label{geomquotchar}
Let $\Theta$ be a bunch in $(E \topto{Q} K, \gamma)$.
The open toric subvariety $X_{\mathfrak{F}(\Theta)} \subset X$ admits a
geometric quotient by the action of $T$ if and only if $\Theta$ is
geometric. 
\end{prop}

\proof $X_{\mathfrak{F}(\Theta)}$  admits a geometric quotient by $T$
if and only if $P \colon F \to N$ is injective on $\vert \Sigma
\vert$. In our situation, the latter is equivalent to saying that
$P \colon F \to N$ is injective on the maximal cones of $\Sigma$,
that means on the cones $\gamma_{0}^{*}$ with $\gamma_{0}$ an element of
$\cov(\Theta)$. 

Thus Lemma~\ref{surjectiveandinjective1} tells us that
$X_{\mathfrak{F}(\Theta)}$ admits a geometric
quotient if and only if every cone $Q(\gamma_{0})$, $\gamma_{0}
\in \cov(\Theta)$, is of full dimension in $K$. Since the elements
of $\Theta$ occur among these cones and for any two cones of
$\cov(\Theta)$ their relative interiors intersect, we obtain the
desired characterization. \endproof

\section{Standard bunches and toric varieties}
\label{section6}

We introduce the class of standard bunches.
The main result of this section, Theorem~\ref{standard2maximal}, 
says that every nondegenerate 2-complete toric variety can be
described by such a standard bunch, and, 
moreover, the isomorphism classes of free standard 
bunches correspond to
the isomorphism classes of nondegenerate 2-complete toric 
varieties having free class group.

\begin{defi}\label{standarddefi}
Let $\Theta$ be a bunch in a projected cone
$(E \topto{Q} K, \gamma)$, and let
$\gamma_{1}, \ldots, \gamma_{n}$ be the facets of 
$\gamma$.
We say that $\Theta$ is a {\em standard bunch\/} if 
\begin{enumerate}
\item for all $i = 1, \ldots, n$ we have $K = Q(\lin(\gamma_{i}) \cap E)$,
\item for every $i = 1, \ldots, n$ there is a $\tau \in \Theta$
  with $\tau^{\circ} \subset Q(\gamma_{i})^{\circ}$.
\end{enumerate}%
If $\Theta$ is a standard bunch in $(E \topto{Q} K, \gamma)$, 
and the cone $\gamma \subset E_{\QQ}$ is regular, then we
speak of the {\em free standard bunch\/} $\Theta$.
\end{defi}

The constructions of Section~\ref{bunchesandadmissiblefans}
provide a functor from standard bunches to toric varieties:

\begin{defi}
Let $\Theta$ be a standard bunch in $(E \topto{Q} K, \gamma)$, 
and let $\Delta := \mathfrak{Q} (\mathfrak{F}(\Theta))$
be the quotient fan of the maximal projectable fan corresponding 
to $\Theta$.
The {\em toric variety associated to\/} $\Theta$ is 
$X_{\Theta} := X_{\Delta}$.
\end{defi}

Recall from Section~\ref{section1} that a toric variety
$X$ is nondegenerate, if it
does not admit a toric decomposition $X \cong X' \times \KK^{*}$.  
Moreover, $X$ is 2-complete, if any toric open embedding 
$X \subset X'$ with $\codim(X' \setminus X) \ge 2$ is an isomorphism.  

\begin{thm}\label{standard2maximal}
The assignment $\mathfrak{T} \colon \Theta \mapsto X_{\Theta}$
defines a contravariant functor from the category of
standard bunches to the category of nondegenerate 2-complete toric
varieties. Moreover,
\begin{enumerate}
\item Every nondegenerate 2-complete toric variety is isomorphic to a
  toric variety $X_{\Theta}$ with a standard bunch $\Theta$.
\item $\mathfrak{T}$ induces a bijection from the isomorphism classes of
  free standard bunches to the isomorphism classes of nondegenerate 
  2-complete toric varieties with free class group.
\end{enumerate}
\end{thm}

For the proof of Theorem~\ref{standard2maximal}, we have to do some
preparation. We need a torsion free version of Cox's
construction~\ref{coxconstr} for nondegenerate fans, 
i.e. fans $\Delta$ in
a lattice $N$ such that the support $\vert \Delta \vert$ generates
the vector space $N_{\QQ}$:

\begin{defi}\label{reducedcox}
Let $(F \topto{P} N, \delta)$ be a projected cone,
$\Sigma$ a fan in $F$, and $\Delta$ a fan in $N$. 
We say that these data form a {\em reduced Cox construction\/}
for $\Delta$ if 
\begin{enumerate}
\item $\Sigma^{(1)}$ equals $\delta^{(1)}$, and $P$ induces 
  bijections $\Sigma^{(1)} \to \Delta^{(1)}$ and 
  $\Sigma^{\max} \to \Delta^{\max}$, 
\item $P$ maps the primitive generators of $\delta$ to primitive
  lattice vectors.
\end{enumerate}
\end{defi}

We show now that every nondegenerate fan admits
reduced Cox contructions. In fact, these reduced
Cox constructions will be even compatible with a 
certain type of maps of fans.

Let $\Delta_{i}$ be nondegenerate fans in lattices $N_{i}$. Moreover,
let $\b{\Psi} \colon N_{1} \to N_{2}$ be any isomorphism of lattices that
is a map of the fans $\Delta_{1}$ and $\Delta_{2}$. Suppose that
$\b{\Psi}$ induces a bijection on the sets of rays $\Delta_{1}^{(1)}$ and
$\Delta_{2}^{(1)}$.

\begin{lemma}\label{lift2cox}
\begin{enumerate}
\item There exist a projected cone
$(F_{1} \topto{P_{1}} N_{1}, \delta_{1})$
and a fan $\Sigma_{1}$ in $F_{1}$ defining a
reduced Cox construction for $\Delta_{1}$.
\item For every reduced Cox construction as in (i)
there exist a reduced Cox construction for 
$\Delta_2$ given by $(F_{2} \topto{P_{2}} N_{2}, \delta_{2})$ 
and $\Sigma_{2}$, 
a lattice isomorphism $\Psi \colon F_{1} \to F_{2}$, 
and a commutative diagram of maps of fans
$$
\xymatrix{
{F_{1}} \ar[r]^{\Psi} \ar[d]_{P_{1}} &
{F_{2}} \ar[d]^{P_{2}} \\
N_{1} \ar[r]_{\b{\Psi}} & N_{2}
}
$$
where the map $\Psi$ induces a bijection of the sets
of rays $\Sigma_{1}^{(1)}$ and $\Sigma_{2}^{(1)}$. 
Moreover, if $\b{\Psi}$ maps $\Delta_{1}$ isomorphically onto a subfan of
$\Delta_{2}$, then $\Psi$ maps $\Sigma_{1}$ isomorphically
onto a subfan of $\Sigma_{2}$.
\end{enumerate}
\end{lemma}

\proof First we perform Cox's original construction~\ref{coxconstr} 
for the fans $\Delta_{i}$. 
Denote by $\Rho_{i}$ the set of rays of $\Delta_{i}$. 
For every maximal cone $\sigma \in \Delta_{i}$ let 
$$ 
\t{\sigma} 
:= \cone(e_{\varrho}; \; \varrho \in \Rho_{i},
                      \; \varrho \subset \sigma). 
$$

Then these cones $\t{\sigma}$ are the maximal cones of a fan
$\t{\Delta}_{i}$ in $\ZZ^{\Rho_{i}}$. Moreover, we have canonical
projections sending the canonical base vectors to the primitive
lattice vectors of the corresponding rays: 
$$ 
C_{i} \colon \ZZ^{\Rho_{i}} \to N_{i},
\qquad
e_{\varrho} \mapsto v_{\varrho}.
$$ 

Let us verify (i). Since $C_{1}$ needs not be surjective, 
it cannot serve as a projection of a projected cone. 
We have to perform a reduction step: 
Let $L_{1} := \ker(C_{1})$, choose a section 
$s_{1} \colon {N_{1}}_{\QQ} \to \QQ^{\Rho_{1}}$ of 
$C_{1}$, and set 
$$ 
F_{1} \; := \; s_{1}(N_{1}) \oplus L_{1} \; \subset \QQ^{\Rho_{1}}.
$$

Then we can view $\delta_{1}$ and $\Sigma_{1}$ as well as objects in
the lattice $F_{1}$. Note that $\delta_{1}$ needs no longer be
regular, but remains simplicial. Together with the surjection
$P_{1} \colon F_{1} \to N_{1}$, the cone $\delta_1$ and the fan 
$\Sigma_1$ give the desired data.

We turn to (ii). Define a (invertible) linear map
$\Psi \colon \QQ^{\Rho_{1}} \to \QQ^{\Rho_{2}}$ 
of rational vector spaces by prescribing its values 
on the canonical base vectors as follows:
$$
\Psi(e_{\varrho}) := e_{\b{\Psi}(\varrho)} \quad
\text{for every } \varrho  \in \Rho_{1}.
$$

Similar to the proof of (i), we may reduce the Cox construction
of $\Delta_2$ by refining the lattice $\ZZ^{\Rho_2}$ via 
$F_{2} :=  \Psi(F_{1})$.
Then the resulting $(F_{2} \topto{P_{2}} N_{2}, \delta_{2})$ 
and $\Sigma_{2}$ and the lattice isomorphism
$\Psi \colon F_1 \to F_2$ are as desired.
\endproof

The above proof shows in particular that we cannot expect uniqueness
of reduced Cox constructions for a given fan.

\begin{lemma}\label{freecl2regular}
Let $\Delta$ be a nondegenerate fan in a lattice $N$. Then the
associated toric variety $X$ has free class group $\Cl(X)$ if and only
if $\Delta$ admits a reduced Cox construction $\Sigma$ in a projected cone
$(F \topto{P} N, \delta)$ with a regular cone $\delta$.
\end{lemma}

\proof Again, we consider Cox's construction~\ref{coxconstr} 
for $\Delta$.
Let $\Rho$ denote the set of rays of $\Delta$, and let $C \colon
\ZZ^{\Rho} \to N$ be the map sending the canonical base vector
$e_{\varrho}$ to the primitive lattice vector $v_{\varrho} \in
\varrho$.

If the toric variety $X$ has free class group then the lattice
homomorphism $C$ is surjective,
see~\cite{Co}. Hence it defines the desired reduced Cox construction
with $F = \ZZ^{\Rho}$ and $\delta$ the positive orthant in
$\QQ^{\Rho}$. 

Conversely, let $(F \topto{P} N, \delta)$ and
$\Sigma$ be a reduced Cox construction of $\Delta$ with $\delta$
regular. Then  $C$ factors as $C = P \circ S$, where $S \colon \ZZ^{\Rho}
\to F$ maps the canonical base vectors to the primitive generators of the
cone $\delta$. It follows that $C$ is surjective. Hence $\Cl(X)$ is
free. \endproof

\begin{lemma}\label{liftisos}
Let $(F_{i} \topto{P_{i}} N_{i}, \delta_{i})$ and $\Sigma_{i}$
be reduced Cox constructions for $\Delta_{i}$ such that the
$\delta_{i}$ are regular cones. Then
every isomorphism $\b{\Psi} \colon N_{1} \to N_{2}$ of the fans
$\Delta_{i}$ admits a unique lifting $\Psi \colon F_{1} \to
F_{2}$ to an isomorphism of the fans $\Sigma_{i}$.
\end{lemma}

\proof For $\varrho \in \Delta_{i}^{(1)}$, let 
$e^{i}_{\varrho} \in F_{i}$ denote the primitive generator of 
$\delta_{i}$
above the primitive vector of $\varrho$. Then define $\Psi \colon
F_{1} \to F_{2}$ by setting $\Psi(e^{1}_{\varrho}) :=
e^{2}_{\b{\Psi}(\varrho)}$. \endproof

Let $\Theta$ be a bunch in $(E \topto{Q} K, \gamma)$,
and let $\Sigma$ denote the associated maximal projectable fan in the
dual projected cone $(F \topto{P} N, \delta)$.
Let $\Delta$ denote the quasifan in $N$ obtained by projecting the
maximal cones of $\Sigma$.

\begin{lemma}\label{standard2cox}
The following statements are equivalent:
\begin{enumerate}
\item $\Delta$ is a fan having $\Sigma$ and
$(F \topto{P} N, \delta)$ as a reduced Cox construction.
\item $\Theta$ is a standard bunch. 
\end{enumerate}
\end{lemma}

\proof
Recall that $\Sigma$ being projectable means that any 
two maximal cones of $\Sigma$ can be separated by a 
$\ker(P)$-invariant linear form. 
In particular, $P$ sets up a bijection 
$\Sigma^{\max} \to \Delta^{\max}$.
Thus, the first statement holds if and only if we have:
\begin{itemize}
\item[(a)] $P$ maps every primitive generator of 
  $\delta$ to a primitive lattice vector in $N$.
\item[(b)] $\Delta$ is a fan,
    we have $\delta^{(1)} = \Sigma^{(1)}$, 
    and $P$ induces a bijection $\Sigma^{(1)} \to \Delta^{(1)}$. 
\end{itemize}

Hence the task is to show that the conditions~\ref{standarddefi}~(i) 
and~(ii) hold if and only if (a) and (b) do so.
Equivalence of (a) and~\ref{standarddefi}~(i) is a direct 
application of Lemma~\ref{surjectiveandinjective2};
note that for this one has to interchange
the roles of $P,\delta$ and $Q,\gamma$  in the lemma such that 
for the $\gamma_{0}$ of the lemma one can take a ray of $\Sigma$.

Now, suppose that the conditions (a) and (b) are valid.
By equivalence of~(a) and~\ref{standarddefi}~(i), 
we only have to check that~\ref{standarddefi}~(ii) 
holds.

Fix a facet $\gamma_{i} \preceq \gamma$.
By assumption, there is a maximal cone 
$\rq{\sigma} \in \Sigma$ such that 
$\rq{\varrho} := \gamma_{i}^{*}$ is a ray of 
$\rq{\sigma}$ and $P(\rq{\varrho})$ is a ray of 
$P(\rq{\sigma}) \in \Delta$.
Then $\gamma_{0} := \rq{\sigma}^{*}$
belongs to $\cov(\Theta)$,
we have $\gamma_{0} \preceq \gamma_{i}$,
and the Invariant Separation Lemma implies
$Q(\gamma_{0})^{\circ} \subset Q(\gamma_{i})^{\circ}$.
By the Overlapping Property~\ref{interior},
any $\tau \in \Theta$ with 
$\tau \subset Q(\gamma_{0})$
is as in~\ref{standarddefi}~(ii).

Conversely, if~\ref{standarddefi}~(i) and~(ii) 
are valid, then we have to show that
$\Delta$ is a fan, $\delta^{(1)}$
equals $\Sigma^{(1)}$ and that $P$ induces 
a bijection $\Sigma^{(1)} \to \Delta^{(1)}$.

Consider any $\rq{\varrho} \in \delta^{(1)}$.
Then $\rq{\varrho} = \gamma_{i}^{*}$ for some facet
$\gamma_{i} \preceq \gamma$.
By~\ref{standarddefi}~(ii), we have
$\tau_{i}^{\circ} \subset Q(\gamma_{i})^{\circ}$ 
for some $\tau_{i} \in \Theta$.
Thus we find a $\gamma_{0} \preceq \gamma_{i}$ 
being minimal with the property that 
$Q(\gamma_{0}) \supset \tau_{0}$ holds 
for some $\tau_{0} \in \Theta$.
Then $\gamma_{0} \in \cov(\Theta)$, and 
$Q(\gamma_{i})^{\circ} \cap Q(\gamma_{0})^{\circ}$
is nonempty because it contains 
$\tau_{i}^{\circ} \cap \tau_{0}^{\circ}$.
Hence $\rq{\sigma} := \gamma_{0}^{*}$
is a maximal cone of $\Sigma$ with
$\rq{\varrho} \preceq \rq{\sigma}$,
and the Invariant Separation Lemma yields 
$P(\rq{\varrho}) \preceq P(\rq{\sigma}) \in \Delta$.

So, this consideration gives in particular
$\delta^{(1)} = \Sigma^{(1)}$.
Moreover, since $P(\rq{\varrho})$ is strictly convex,
it gives $0 \in \Delta$; in other words,
the quasifan $\Delta$ is a fan.
Furthermore, since we already know that (a) holds,
the image $P(\rq{\varrho})$ is in fact one-dimensional. 
Hence we obtain $P(\rq{\varrho}) \in \Delta^{(1)}$,
and thus $\Sigma^{(1)} \to \Delta^{(1)}$
is well defined.

Surjectivity of the map $\Sigma^{(1)} \to \Delta^{(1)}$
follows from the fact that 
$\Sigma^{\max} \to \Delta^{\max}$ is surjective.
Injectivity follows from the observation that 
by~\ref{standarddefi}~(ii) we always have
$Q(\gamma_{i})^{\circ} \cap Q(\gamma_{j})^{\circ} \ne \emptyset$, 
and hence any two rays of $\delta$ can be 
invariantly separated.
\endproof

\proof[Proof of Theorem~\ref{standard2maximal}]
By Lemma~\ref{standard2cox}, the toric variety $X_\Theta$ 
associated to a standard bunch $\Theta$ is nondegenerate.
It is also 2-complete: Otherwise there is an open toric embedding
$X_{\Theta} \subset X$ with nonempty complement of
codimension at least two. 
Using Lemmas~\ref{lift2cox} and~\ref{standard2cox}, 
we can compare reduced Cox constructions of
$X_{\Theta}$ and $X$, and we see that 
the projectable fan associated to $\Theta$ 
is not maximal, i.e., does not
satisfy~\ref{admissiblefandefi}~(i)~(b).
A contradiction.
 
So the functor $\mathfrak{T} \colon \Theta \mapsto X_{\Theta}$ is
well defined. The fact that it is surjective on isomorphism classes
follows from existence of reduced Cox constructions,
Proposition~\ref{goodquotchar}, Theorem~\ref{dualitythm} and
Lemma~\ref{standard2cox}. The correspondence
of isomorphism classes of free bunches with isomorphism classes of
nondegenerate 2-complete toric varieties with free class group is a direct
application of Lemmas~\ref{freecl2regular} and~\ref{liftisos}. 
\endproof

\section{A very first dictionary}\label{section7}

Fix a standard bunch $\Theta$ in a projected cone 
$(E \topto{Q} K, \gamma)$, 
and let $X := X_{\Theta}$ denote the associated toric
variety.  
In this section, we characterize basic geometric
properties of $X$ in terms of the bunch $\Theta$. 

Let $(F \topto{P} N, \delta)$ be the dual projected cone. Denote by $\Sigma$
the maximal projectable fan associated to $\Theta$, and let $\Delta$ be the
quotient fan of $\Sigma$. Recall from Lemma~\ref{standard2cox} that these data
form a reduced Cox construction of $\Delta$. In particular, $\Delta$ lives in
$N$, and we have 
$$X = X_{\Delta}, \qquad \dim (X) = \rank (E) - \rank (K).$$

We study now $\QQ$-factoriality, smoothness, existence of fixed points and
completeness. For this we need the following observation:

\begin{lemma}\label{fullandsimplicial}
Consider a face $\gamma_{0} \in \cov(\Theta)$ and the corresponding
maximal cone $P(\gamma_{0}^{*})$ of $\Delta$. Then we have:
\begin{enumerate}
\item $Q(\gamma_{0})$ is of full dimension if and only if
  $P(\gamma_{0}^{*})$ is simplicial.
\item $Q(\gamma_{0})$ is simplicial if and only if
  $P(\gamma_{0}^{*})$ is of full dimension.
\end{enumerate}
\end{lemma}

\proof We prove (i). Let $Q(\gamma_{0})$ be of full dimension. By
Lemma~\ref{surjectiveandinjective1}, the map $P$ is injective
along $\gamma_{0}^{*}$. In particular, $P(\gamma_{0}^{*})$ is
simplicial. Conversely, let $P(\gamma_{0}^{*})$ be simplicial. Since
$P$ induces a bijection from the rays of $\gamma_{0}^{*}$ to
the rays of $P(\gamma_{0}^{*})$, it is injective along
$\gamma_{0}^{*}$. Thus Lemma~\ref{surjectiveandinjective1} yields that 
$Q(\gamma_{0})$ is of full dimension.

We turn to (ii). If $P(\gamma_{0}^{*})$ is of full dimension, we see
as before that $Q(\gamma_{0})$ is simplicial. For the converse we show
that $Q$ is injective along $\gamma_{0}$: For every ray $\varrho$ of
$Q(\gamma_{0})$, choose a ray
$\tau$ of $\gamma_{0}$ with $Q(\tau) = \varrho$. Then the cone
$\gamma_{1} \preceq \gamma_{0}$ generated by these rays $\tau$ is mapped
bijectively onto $Q(\gamma_{0})$. By minimality of $\gamma_{0}$ as an
element of $\cov(\Theta)$, we conclude $\gamma_{1} = \gamma_{0}$.
\endproof

The first statement of this section is the following
characterization of $\QQ$-facto\-riality:

\begin{prop}\label{Qfactorialchar}
The toric variety $X$ is $\QQ$-factorial if and only if
$\Theta$ consists of cones of full dimension in $K$.
\end{prop}

\proof The toric variety $X$ is $\QQ$-factorial if and only if
all cones of $\Delta$ are simplicial. By Lemma~\ref{fullandsimplicial}
this is equivalent to saying that all cones $Q(\gamma_{0})$,
$\gamma_{0} \in \cov(\Theta)$, are of full dimension. 
The latter holds if and only if $\Theta$ consists of cones of 
full dimension, because every $Q(\gamma_{0})$ contains some 
$\tau \in \Theta$.
\endproof

Characterizing smoothness of  toric varieties in terms of bunches
involves the lattice structure:
  
\begin{prop}\label{smoothchar}
The toric variety $X_{\Theta}$ is smooth if and only if for every 
$\gamma_{0} \in \cov(\Theta)$ we have
\begin{enumerate}
\item $\gamma_{0}^{*}$ is a regular cone in $F$.
\item $Q$ maps $\lin(\gamma_{0}) \cap E$ onto $K$.
\end{enumerate}
\end{prop}

\proof Suppose that (i) and (ii) hold. To verify smoothness of $X$, we
have to show that all cones  $P(\gamma_{0}^{*})$, $\gamma_{0} \in
\cov(\Theta)$, are regular. By the properties~\ref{reducedcox} 
of a reduced Cox construction and Condition~(i),
the primitive generators of $P(\gamma_{0}^{*})$ span the
sublattice $P(\lin(\gamma_{0}^{*}) \cap F)$ of $N$.
By Lemma~\ref{surjectiveandinjective2} and Condition~(ii), this
sublattice is primitive. This proves regularity of the cone
$P(\gamma_{0}^{*})$. 

Conversely, let $X$ be smooth. For $\gamma_{0} \in
\cov(\Theta)$, consider the sublattice $F_{0} \subset F$ spanned by
the primitive generators of $\gamma_{0}^{*}$. Since $P(\gamma_{0}^{*})$
is regular, the properties~\ref{reducedcox} of a reduced Cox
construction yield that $P(F_{0})$ is a 
primitive sublattice of $N$. Since $P$ is injective along
$\gamma_{0}^{*}$, we see that also $F_{0}$ is primitive. This gives
Condition~(i). Condition~(ii) then follows from 
Lemma~\ref{surjectiveandinjective2}. \endproof

Existence of global regular functions on a toric variety is
characterized as follows:

\begin{prop}\label{nontrivfunctionschar} 
We have $\mathcal{O}(X) = \KK$ if and only if $Q$ contracts no ray of
$\gamma$ to a point and the image $Q(\gamma)$ is strictly convex. 
\end{prop}

\proof 
$\mathcal{O}(X) = \KK$ holds if and only if the rays of 
$\Delta$ generate $N_{\QQ}$ as a convex cone. 
This is valid if and only if $\delta + \ker(P)_{\QQ}$ 
equals $F_{\QQ}$. 
By dualizing, this condition is equivalent to
$\gamma \cap \ker(Q)_{\QQ} = \{0\}$. 
This in turn holds if and only if $Q$ contracts
no ray of $\gamma$ and $Q(\gamma)$ is strictly convex.
\endproof

Recall that we speak of a full toric variety $X$ if $X$ is
2-complete and every $T_{X}$-orbit contains a fixed point in its closure.

\begin{prop}\label{existenoffixedpointschar}
The toric variety $X$ is full if and only if all cones
$Q(\gamma_{0})$, where $\gamma_{0} \in \cov(\Theta)$, 
are simplicial.
\end{prop}

\proof 
Existence of fixed points in the $T_{X}$-orbit closures means 
that all maximal cones of $\Delta$ are of full
dimension. Thus Lemma~\ref{fullandsimplicial} gives the assertion.
\endproof

Finally, we also can characterize completeness of a toric variety in
terms of its bunch:

\begin{prop}\label{completenesschar}
The toric variety $X_{\Theta}$ is complete if and only if $\Theta$
contains a simplicial cone and any face $\gamma_{0} \preceq \gamma$
satisfying $\tau^{\circ} \subset Q(\gamma_{0})^{\circ}$ 
for some $\tau \in \Theta$ and $\gamma_{1} \preceq \gamma_{0}$ 
for only one $\gamma_{1} \in \cov(\Theta)$ belongs to 
$\cov(\Theta)$.  
\end{prop}

\proof
This is a direct translation of a well known characterization of
completeness in terms of fans. Namely, $X$ is complete if and
only the fan $\Delta$ has the following two properties: firstly,
at least one of its (maximal) cones is of full dimension; secondly, 
any cone contained in only one maximal cone is itself
maximal.

By Lemma~\ref{fullandsimplicial}, the first property translates 
to the property that $\Theta$ contains a simplicial cone.
For the second, recall from Proposition~\ref{conesofquotfan}, 
that the cones $\sigma_{0} \in \Delta$ correspond to 
the faces $\gamma_{0} \preceq \gamma$ satisfying 
$\tau^{\circ} \subset Q(\gamma_{0})^{\circ}$ for some 
$\tau \in \Theta$ via
$$
\sigma_{0} \mapsto 
\left(P^{-1}(\sigma_{0}) \cap \vert \Sigma \vert\right)^{*}.
$$ 
Thereby the maximal cones $\sigma_{1}$ of $\Delta$ correspond to
the elements of $\cov(\Theta)$. Thus the statement that 
$\sigma_{0} \preceq \sigma_{1}$ for only one maximal 
$\sigma_{1}$ implies $\sigma_{0} \in \Delta^{\max}$ 
directly translates to the second characterizing condition 
of the assertion.
\endproof

One may ask if a full toric variety $X$ with 
$\mathcal{O}(X) = \KK$ is already complete. 
If $\dim(X) \le 3$ holds, then the answer is positive. 
For $\dim(X) \ge 4$, there are counterexamples,
compare~\cite[Remark~2]{Ew1}.
However, we shall see in Section~\ref{section10}
that every smooth 2-complete toric variety
with class group $\ZZ^{2}$ is complete.

\section{Full $\QQ$-factorial toric varieties}\label{section8}

In this section we introduce the class of simple bunches, 
and we show that $\Theta \mapsto X_{\Theta}$ defines an equivalence 
of categories between the simple bunches and the 
full $\QQ$-factorial toric varieties.

\begin{defi}
By a {\em simple bunch\/} we mean a standard bunch $\Theta$ in a
projected cone $(E \topto{Q} K, \gamma)$ such that
$Q$ maps $E \cap \lin(\gamma_{0})$ isomorphically onto $K$ for every
$\gamma_{0} \in \cov(\Theta)$.
\end{defi}

Note that the cones of a simple bunch are of full dimension
and simplicial, but they need not be regular. 
To state the main result of this section,
recall from Definition~\ref{nondegfulldef}
that a toric variety $X$ is full, if it is
2-complete and every $T_{X}$-orbit has a fixed point in its closure.

\begin{thm}\label{equivcat}
The assignment $\Theta \mapsto X_{\Theta}$ defines an equivalence of
the category of simple bunches with the category of full 
$\QQ$-factorial toric varieties.
\end{thm}

The crucial observation for the proof is that full $\QQ$-factorial
toric varieties admit a universal reduced 
Cox construction in a certain sense. 
First we show existence:
 
\begin{lemma}\label{freeexists}
Let $\Delta$ be a simplicial fan in a lattice $N$ such that any maximal cone
is of full dimension. Then there are $(F \topto{P} N, \delta)$ and
$\Sigma$ defining a reduced Cox construction of $\Delta$ such that 
\begin{enumerate}
\item $P$ induces for every $\rq{\sigma} \in \Sigma^{\max}$ an isomorphism
$P_{\rq{\sigma}} \colon F \cap \lin(\rq{\sigma}) \to N$.
\item $F$ is the sum of the sublattices $\lin(\rq{\sigma}) \cap F$,
where $\rq{\sigma} \in \Sigma^{\max}$.
\end{enumerate} 
\end{lemma}
 
\proof As usual, let $\Rho$ be the set of rays of $\Delta$, let 
$P \colon \ZZ^{\Rho} \to N$ be the homomorphism sending the canonical
base vector $e_{\varrho}$ to the primitive vector $v_{\varrho} \in
\varrho$, and let $\delta \subset \QQ^{\Rho}$ 
be the cone generated by the $e_{\varrho}$.
For any cone $\sigma \in \Delta$ set
$$ 
\rq{\sigma} 
\; := \; 
\cone(e_{\varrho}; \; \varrho \in \sigma^{(1)}). 
$$

Then these cones form a fan $\Sigma$ in the lattice $\ZZ^{\Rho}$. 
In order to achieve the desired properties, we refine the lattice 
$\ZZ^{\Rho}$ as follows: 
Think of $P$ for the moment as a map of the vector spaces
$\QQ^{\Rho}$ and $N_{\QQ}$ and consider the set
$$ 
P^{-1}(N) \cap \vert \Sigma \vert
\; = \; \bigcup_{\sigma \in \Sigma} P^{-1}(N) \cap \rq{\sigma}
\; \subset \;
\QQ^{\Rho}. 
$$

This set generates a lattice $F \subset \QQ^{\Rho}$,
because by simpliciality of $\Delta$ the map $P$ is injective
along the cones $\rq{\sigma} \in \Sigma$, and hence each
$\rq{\sigma} \cap P^{-1}(N)$ is discrete.
Moreover, the restriction $P_{\rq{\sigma}}
\colon F \cap \lin(\rq{\sigma}) \to N$ of $P$ is an isomorphism 
for every $\rq{\sigma} \in \Sigma^{\max}$. 
Here we use that the maximal cones of $\Delta$ are of full 
dimension.

We may view $\delta$ and $\Sigma$ as well as data in the lattice $F$. 
Since $P \colon F \to N$ is surjective, this gives in particular
a projected cone $(F \topto{P} N, \delta)$. Now it is 
straightforward to verify the defining properties
of a reduced Cox construction for these data.
By construction it satisfies~(i) and~(ii). 
\endproof 

We shall call a reduced Cox construction with the properties
of Lemma~\ref{freeexists} a 
{\em universal reduced Cox construction}. 
Uniqueness of universal reduced Cox constructions is a 
consequence of the following lifting property: 

\begin{lemma}\label{liftinglemma}
Let $\Delta_{i}$ be fans in lattices $N_{i}$, and let $(F_{i} \topto{P_{i}}
N_{i}, \delta_{i})$ and $\Sigma_{i}$ be universal reduced Cox
constructions. Then every map $\b{\Psi} \colon N_{1} \to N_{2}$ of the fans
$\Delta_{i}$ admits a unique lifting to a map $\Psi \colon F_{1} \to F_{2}$ of
the fans $\Sigma_{i}$:
$$ \xymatrix{ {F_{1}} \ar[r]^{\Psi} \ar[d]_{{P_{1}}} & F_{2}
\ar[d]^{{P_{2}}} \\ {N_{1}} \ar[r]_{\b{\Psi}} & N_{2} }$$
\end{lemma}

\proof 
For every maximal cone $\sigma_{1} \in \Delta_{1}$, 
fix a maximal cone $\sigma_{2} \in \Delta_{2}$ 
with $\b{\Psi}(\sigma_{1}) \subset \sigma_{2}$. 
Let $\rq{\sigma}_{i} \in \Sigma_{i}$ denote the cones lying over
$\sigma_{i}$. 
Then, by the Property~\ref{freeexists}~(i) of a 
universal reduced Cox construction,
we obtain for every $\sigma_{1}$ a unique commutative diagram
$$ 
\xymatrix{ 
{F_{1} \cap \lin(\rq{\sigma}_{1})}
\ar[r]^{\Psi_{\rq{\sigma}_{1}}} \ar[d]_{{P_{1}}} 
& 
F_{2} \cap \lin(\rq{\sigma}_{2}) \ar[d]^{{P_{2}}} \\ 
{N_{1}} \ar[u]_{s_{\rq{\sigma}_{1}}}
\ar[r]_{\b{\Psi}} 
& 
N_{2} \ar[u]^{s_{\rq{\sigma}_{2}}} }
$$ 
where the $s_{\rq{\sigma}_{i}}$ are the sections
mapping the primitive generators of $\sigma_{i}$ to those 
of~$\rq{\sigma}_{i}$.
Any two $\Psi_{\rq{\sigma}_{1}}$ and
$\Psi_{\rq{\sigma}'_{1}}$ have the same values on
the primitive generators of $\delta$ that lie in
$\rq{\sigma}_{1} \cap \rq{\sigma}'_{1}$.
Thus the $\Psi_{\rq{\sigma}_{1}}$ fit together to
a linear map $\Psi \colon F_{1} \to (F_{2})_{\QQ}$.
Lemma~\ref{freeexists}~(ii) gives
$\Psi(F_{1}) \subset F_{2}$.
Hence we found the desired lifting.
\endproof 

\proof[Proof of Theorem \ref{equivcat}] 
First observe that for a simple bunch $\Theta$, the toric 
variety $X_{\Theta}$ is indeed full and $\QQ$-factorial: 
Since for every $\gamma_{0} \in \cov(\Theta)$ the cone $Q(\gamma_{0})$ 
is simplicial and of full dimension, Lemma~\ref{fullandsimplicial} 
ensures that also the maximal cones of the fan $\Delta$ defining $X$ 
are of full dimension and simplicial.

Next we show that, up to isomorphism, every full $\QQ$-factorial toric
variety $X$ is of the form $X_{\Theta}$ with a simple bunch. We may
assume that $X$ arises from a fan $\Delta$ in a lattice $N$. Lemma
\ref{freeexists} provides a universal reduced Cox construction 
$(F \topto{P} N, \delta)$ and $\Sigma$ of $\Delta$. By
Lemma~\ref{surjectiveandinjective2}, the bunch $\Theta$ corresponding
to the maximal projectable fan $\Sigma$ is as wanted.

Finally, we have to show bijectivity on the level of morphisms. Due to
Lemma~\ref{liftinglemma} the morphisms between full $\QQ$-factorial 
simplicial toric varieties are in one-to-one correspondence with the
maps of maximal projectable fans defined by their 
universal reduced Cox constructions.
By Theorem~\ref{dualitythm}, the latter morphisms correspond to the
morphisms of the associated bunches.   
\endproof

\begin{coro}
$\Theta \mapsto X_{\Theta}$ defines an equivalence of categories from
free simple bunches to full smooth toric varieties.
\end{coro}

\proof
By Propositions~\ref{smoothchar} and~\ref{existenoffixedpointschar}, 
every free simple bunch $\Theta$ defines a full smooth toric variety
$X_{\Theta}$.
Conversely, for a full smooth toric variety $X$, 
the usual Cox construction~\ref{coxconstr} is a universal 
reduced Cox construction.
Hence $X$ arises from a free simple bunch.
\endproof

We conclude this section with some observations and examples. 
The first remark gives the geometric interpretation of the 
universal reduced Cox construction: 

\begin{rem}
In terms of toric varieties, the universal reduced Cox
construction $p \colon \rq{X} \to X$ arises from the usual 
Cox construction $c \colon \t{X} \to X$ by dividing $\t{X}$ 
by the (finite) group generated by all
isotropy groups of the action of $\ker(Q)$ on $\t{X}$.
\end{rem}

The following example shows that for nonsimplicial toric varieties
the functor $\Theta \mapsto X_{\Theta}$ associating to a 
standard bunch its toric variety is not surjective on the level of
morphisms:

\begin{exam}\label{freegroup}
We present a toric morphism that cannot be lifted to the Cox
constructions. In $\ZZ^{3}$ consider the vectors
$$ 
v_{1} := (1,0,0),
\quad
v_{2} := (1,0,1),
\quad
v_{3} := (0,1,0),
\quad
v_{4} := (-3,2,2),
\quad
v_{5} := (0,1,1).
$$

Let $\Delta_{2}$ be the fan of faces of the cone generated by $v_{1},
\ldots, v_{4}$, and let $\Delta_{1}$ be the subdivision of
$\Delta_{2}$ at $v_{5}$. Mapping the $i$-th canonical base vector to
$v_{i}$, we obtain projections
$$ 
P_{1} \colon \ZZ^{5} \to \ZZ^{3},
\qquad
P_{2} \colon \ZZ^{4} \to \ZZ^{3}.
$$
Let $\Sigma_{1}$ and $\Sigma_{2}$ be the fans above
$\Delta_{1}$ and $\Delta_{2}$, respectively. Note that these data are
in fact reduced Cox constructions.

We claim that the identity $\varphi \colon \ZZ^{3} \to \ZZ^{3}$ does
not admit a lifting. In fact, since $v_{5}$ equals $3/2 v_{1} +
1/2 v_{4}$, a possible lifting $\Phi \colon \ZZ^{5} \to
\ZZ^{4}$ must satisfy
$$ 
\Phi(0,0,0,0,1) 
\; \in \; 
(3/2, 0,0, 1/2) + \QQ(5,-2,-2,1),
$$ 
where $(5,-2,-2,1)$ generates $\ker(P_{2})$. An explicit calculation
shows that the right hand sice does not contain integral points with
nonnegative coefficients. 

This excludes existence of a map $\Phi \colon \ZZ^{5} \to \ZZ^{4}$ of
the fans $\Sigma_{1}$ and $\Sigma_{2}$ lifting $\varphi
\colon \ZZ^{3} \to \ZZ^{3}$.
\end{exam}

The next example shows that in the nonsimplicial case the functor
$\Theta \mapsto X_{\Theta}$
is not injective on the level of morphisms.

\begin{exam}\label{notfaithfull}
We give a toric morphism that admits two different liftings.
Let $\Delta_{2}$ be the fan of faces of the cone generated by
the vectors
$$
(1,0,0),
\quad
(0,1,0),
\quad
(-1,0,1),
\quad
(0,-1,1).
$$

Let $\Delta_{1}$ be the subdivision of $\Delta_{2}$ at $(0,0,1)$. 
Then there are two different liftings of the identity map $\ZZ^{3} \to
\ZZ^{3}$ to the respective Cox-constructions, namely the maps $\ZZ^{5}
\to \ZZ^{4}$ defined by the matrices 
$$ 
\begin{pmatrix}
  1 & 0 & 0 & 0 & 0 \\ 0 & 1 & 0 & 0 & 1 \\ 0 & 0 & 1 & 0 & 0 \\ 0 & 0 & 0 & 1
  & 1 
\end{pmatrix},
\qquad 
\begin{pmatrix} 1 & 0 & 0 &
  0 & 1 \\ 0 & 1 & 0 & 0 & 0 \\ 0 & 0 & 1 & 0 & 1 \\ 0 & 0 & 0 & 1 & 0
  \end{pmatrix}.
$$ 
\end{exam}

\section{Invariant divisors and divisor classes}\label{section9}

In this section, we come to one of the most powerful parts
of the language of bunches: We study geometric properties of
divisors and divisor classes. 

Let us fix the notation. As usual, $\Theta$ is a standard bunch
in a projected cone $(E \topto{Q} K, \gamma)$, 
and $X := X_{\Theta}$ is the associated toric variety
with its big torus $T := T_{X}$.
Let $(F \topto{P} N, \delta)$ be the dual projected cone, 
$\Sigma$ the maximal projectable fan corresponding to $\Theta$, 
and $\Delta$ the quotient fan of $\Sigma$. 
Recall that these data define a reduced Cox construction
of $\Delta$, and we have $X = X_{\Delta}$.

Our first task is to relate the lattice $K$ to the divisor class 
group $\ClDiv(X)$.
Let $v_{1}, \dots, v_{r}$ be the primitive generators of the
one-dimensional cones of $\Delta$, 
and let $\rq{v}_{1}, \dots, \rq{v}_{r}$ 
be the primitive generators of the rays of $\delta$,
numbered in such a way that we always have 
$P(\rq{v}_{i}) = v_{i}$. 

Every ray $\varrho_{i} = \QQ_{\ge 0} v_{i}$ determines an invariant prime
divisor $D_{i}$ in $X$.
There is a canonical injection mapping $E$ into the lattice 
$\WDiv^{T}(X)$ of invariant Weil divisors on $X$:
$$ 
\mathfrak{D} \colon E \to \WDiv^{T}(X), 
\qquad 
\rq{w} \; \mapsto \; 
\mathfrak{D}(\rq{w}) 
:= 
\sum_{i=1}^{r} \rq{w} (\rq{v}_{i}) D_{i}.
$$
By construction, an element $u \in M$ is mapped to the principal
divisor $\div (\chi^{u})$ of $X$. Hence, we obtain, compare~\cite{Fu}:

\begin{prop}\label{divisors}
There is a commutative diagram with exact rows and injective upwards 
arrows:
$$ 
\xymatrix{ 0 \ar[r] & \PDiv^{T}(X) \ar[r] & {\WDiv^{T} (X)} \ar[r] &
\ClDiv(X) \ar[r] & 0 \\ 0 \ar[r] & M \ar[r] \ar[u]^{\cong} & E \ar[r]
\ar[u]^{\mathfrak{D}} & K \ar[r] \ar[u] ^{\b{\mathfrak{D}}} & 0 } 
$$
Moreover, tensoring this diagram with $\QQ$ provides
isomorphisms of rational vector spaces:
$$ 
\mathfrak{D} \colon E_{\QQ} \to \WDiv^{T}_{\QQ} (X),
\qquad
\b{\mathfrak{D}} \colon K_{\QQ} \to \ClDiv_{\QQ}(X).
$$
\end{prop}

In the main result of this section, we study
the group $\Pic_{\QQ}(X) \subset \ClDiv_{\QQ}(X)$ 
of rational Cartier divisor classes.
We obtain very simple descriptions 
of the cone $\C^{\rm sa}(X) \subset \Pic_{\QQ}(X)$ 
of semiample classes
and the cone $\C^{\rm a}(X) \subset \Pic_{\QQ}(X)$ 
of ample classes:

\begin{thm}\label{QCartier}
The map $\b{\mathfrak{D}} \colon K_{\QQ} \to \ClDiv_{\QQ}(X)$
defines canonical isomorphisms:
$$
\Pic_{\QQ}(X) \cong \bigcap_{\tau \in \Theta} \lin(\tau),
 \qquad
\C^{\rm sa}(X) \cong \bigcap_{\tau \in \Theta} \tau,
 \qquad
\C^{\rm a}(X) \cong \bigcap_{\tau \in \Theta} \tau^{\circ}.
$$
\end{thm}

\proof
Let $D \in \WDiv_{\QQ}^{T}(X)$, 
and set $\rq{w} := \mathfrak{D}^{-1}(D) \in E_{\QQ}$.
Our task is to characterize the statements that $D$ 
is a $\QQ$-Cartier, a semiample or 
an ample divisor 
in terms of the image $w :=  Q(\rq{w})$ in $K_{\QQ}$.

For the description of $\Pic(X)$,
recall from~\cite[p.~66]{Fu}, that $D$ is $\QQ$-Cartier
if and only if it arises from a support function, 
i.e., there is a family $(u_{\sigma})_{\sigma \in \Delta^{\max}}$ 
with $u_{\sigma} \in M_{\QQ}$ such that 
$D = m^{-1}\div(\chi^{mu_\sigma})$ holds on each affine chart
$X_\sigma \subset X$ with a postive integral multiple $m u_{\sigma}$.
If $D$ is $\QQ$-Cartier, then the describing
linear forms $u_\sigma$ are unique up to elements of 
$\sigma^\perp$.

Now, suppose that $D$ is $\QQ$-Cartier, and 
let $(u_\sigma)$ be a describing support function. 
Then $u_{\sigma} (v_{i}) = \rq{w}(\rq{v}_{i})$ holds
for every $v_{i} \in \sigma^{(1)}$.
Define $\ell_{\sigma} := \rq{w} - P^{*}(u_{\sigma})$ 
and denote the cone above
$\sigma$ by $\rq{\sigma}$.
Then $\ell_{\sigma}$ lies in 
$\rq{\sigma}^{\perp} = \lin(\rq{\sigma}^{*})$, 
and hence 
$Q(\rq{w}) = Q(\ell_{\sigma})$
lies in $Q(\lin(\rq{\sigma}^{*}))$.
Since this applies for all $\sigma \in \Delta^{\max}$, 
we obtain 
$$
w
\; = \; 
Q(\rq{w})
\; \in \;
\bigcap_{\gamma_{0} \in \cov(\Theta)} \lin(Q(\gamma_{0})) 
\; = \;
\bigcap_{\tau \in \Theta} \lin(\tau).
$$

Conversely, let $w = Q(\rq{w})$ 
belong to the last intersection.
Since for $\sigma \in \Delta^{\max}$ the image 
$Q(\rq{\sigma}^{*})$ contains an element of
$\Theta$, 
we find for each $\sigma \in \Delta^{\max}$ an
$\ell_{\sigma} \in \lin(\rq{\sigma}^{*})$
with $Q(\ell_{\sigma}) = w$. 
Then $u_{\sigma} := \rq{w} - \ell_{\sigma}$ 
maps to zero and can therefore be viewed as an element of
$M_{\QQ}$. 
Hence, $(u_{\sigma})_{\sigma \in \Delta^{\max}}$ provides 
a support function describing $D$.

For the descriptions of $\C^{\rm sa}(X)$ and 
$\C^{\rm a}(X)$, let $D$ be $\QQ$-Cartier.
Recall that $D$ is semiample (ample) if and only if 
it is described by a support function $(u_{\sigma})$, 
which is convex (strictly convex) in the sense that
$u_{\sigma} - u_{\sigma'}$ is nonnegative (positive) 
on $\sigma \setminus \sigma'$ for any two
$\sigma, \sigma' \in \Delta^{\max}$. 

Suppose that $D$ is semiample (ample) with convex 
(strictly convex) support function $(u_{\sigma})$.
In terms of $\ell_{\sigma} := \rq{w} - P^{*}(u_{\sigma})$
this means that each $\ell_{\sigma'} - \ell_{\sigma}$ 
is nonnegative (positive) 
on $\rq{\sigma} \setminus \rq{\sigma}'$.
Since $\ell_{\sigma} \in \rq{\sigma}^{\perp}$
holds, this is equivalent to nonnegativity (positivity) 
of $\ell_{\sigma'}$ on every $\rq{\sigma} \setminus \rq{\sigma}'$.

Since all rays of the cone $\delta$ occur in the fan $\Sigma$,
the latter is valid if and only if
$\ell_{\sigma} \in \rq{\sigma}^{*}$ 
(resp. $\ell_{\sigma} \in (\rq{\sigma}^{*})^{\circ}$)
holds for all $\sigma$. This in turn implies
that for every $\sigma \in \Delta^{\max}$ we have
\begin{equation}
\label{intheimages}
w 
= 
Q(\rq{w}) 
= 
Q(\ell_{\sigma}) 
\in Q(\rq{\sigma}^{*}) 
\quad 
\text{(resp. $w\in Q((\rq{\sigma}^{*})^{\circ})$)}.
\end{equation}

Now, the $\rq{\sigma}^{*}$, where $\sigma \in \Delta^{\max}$, 
are precisely the cones of $\cov(\Theta)$. 
Since any interior $Q(\rq{\sigma}^{*})$ contains the interior 
of a cone of $\Theta$, we can conclude that $w$ lies in the 
respective intersections of the assertion.

Conversely, if $w$ belongs to one of the right hand side 
intersections, then we surely arrive at~\ref{intheimages}.
Thus, for every $\sigma \in \Delta^{\max}$,
we find an $\ell_\sigma \in \rq{\sigma}^*$ 
(an $\ell_{\sigma} \in (\rq{\sigma}^{*})^{\circ}$)
mapping to $w$. 
Reversing the above arguments, 
we see that $u_{\sigma} := \rq{w} - \ell_{\sigma}$
is a convex (strictly convex) support function describing $D$.
\endproof

As an immediate consequence of Theorem~\ref{QCartier}, 
we obtain a quasiprojectivity criterion
in the spirit of~\cite{Ew}.

\begin{coro}\label{Ewaldcrit}
The variety $X$ is quasiprojective if and only if the intersection
of all relative interiors $\tau^{\circ}$, where $\tau \in \Theta$,
is nonempty. \endproof
\end{coro}

Combining Theorem~\ref{QCartier} with Proposition~\ref{Qfactorialchar} 
shows that a $\QQ$-factorial quasiprojective toric variety has an
ample cone of full dimension. Hence we get back a 
result of Oda and Park~\cite[Theorem~3.5]{OP}:

\begin{coro}
Suppose that $X$ is quasiprojective and $\QQ$-factorial.
Then every rational Weil divisor $D$ of $X$ admits a 
representation $D = D_{1} - D_{2}$ with ample rational
Cartier divisors $D_{1}$ and $D_{2}$. \endproof
\end{coro}

So far, we considered rational divisors and divisor classes.
For toric varieties with free class group, 
we also obtain a simple picture for integral divisors.
We need the following notation:

Every ray $\rq{\varrho}_{i} := \QQ_{\ge 0} \rq{v}_{i}$ 
of $\delta$ has a unique ``opposite'' ray, namely
the unique ray of $\gamma$ that is not contained 
in $\rq{\varrho}_{i}^{\bot}$. 
We denote the primitive generator of this opposite ray by
$\rq{w}_{i}$, and its image in $K$ by $w_{i} := Q(\rq{w}_{i})$.

\begin{rem}\label{dualbase}
The bunch $\Theta$ is free if and only if 
$\rq{w}_{i} (\rq{v}_{i}) = 1$ holds
for all $1 \le i \le r$.
\end{rem}

In the free case, we have the following integral version of the
corresponding statements in Proposition~\ref{divisors} and 
Theorem~\ref{QCartier}:

\begin{prop}
Assume that $\Theta$ is a free bunch. Then we have:
$$
K \; \cong \; \ClDiv(X),
\qquad
\Pic(X) 
\; \cong \; 
\bigcap_{\gamma_{0} \in \cov(\Theta)} Q(\lin(\gamma_{0}) \cap E).
$$
\end{prop}

\proof
Using Remark~\ref{dualbase}, we infer from freeness of $\Theta$
that the map $\mathfrak{D} \colon E \to \WDiv^{T}(X)$ is an
isomorphism. 
Using the 5-Lemma in the diagram of Proposition~\ref{divisors}, 
we obtain that also the map
$\b{\mathfrak{D}} \colon K \to \ClDiv(X)$ 
is an isomorphism. 
This gives the first part of the assertion.

The rest of the proof is similar to that of Proposition \ref{QCartier}:
Suppose that $w \in K$ lies in all $Q(\lin(\gamma_{0}) \cap E)$, 
where $\gamma_{0} \in \cov(\Theta)$.
Let $\rq{w} \in E$ with $Q(\rq{w}) = w$.
Then for every $\sigma \in \Delta^{\max}$ we can choose
$\ell_{\sigma}$ even in $\rq{\sigma}^{\perp} \cap E$ with 
$Q(\ell_{\sigma}) = w$.
But then all $u_{\sigma} := \rq{w} - \ell_{\sigma}$ lie
in $M$.
Therefore $\mathfrak{D}(\rq{w})$ is a Cartier divisor and 
thus we have $\b{\mathfrak{D}}(w) \in \Pic(X)$. 

Conversely, given $w \in K$ with $\b{\mathfrak{D}}(w) \in \Pic(X)$,
choose $\rq{w} \in E$ with $Q(\rq{w}) = w$.  
Then $\mathfrak{D}(\rq{w})$ is a Cartier divisor, and hence
it is described by a support function 
$(u_{\sigma})_{\sigma \in \Delta^{\max}}$ with $u_{\sigma} \in M$. 
But then each $\ell_{\sigma} := \rq{w} - P^{*} (u_{\sigma})$ 
lies in $\rq{\sigma}^\perp \cap E$, which proves the assertion.
\endproof

We present some further applications.
If $\Theta$ is free, then we can easily
describe the canonical divisor class
of $X$: By~\cite[Sec.~4.3]{Fu}, the negative 
of the sum  over all invariant prime divisors 
is a canonical divisor on $X$; its class is 
given by
$$ 
[K_X] = - \b{\mathfrak{D}}(w_1 + \ldots + w_r).
$$

Recall that a toric variety is said to be {\em $\QQ$-Gorenstein\/} if
some multiple of its anticanonical divisor is Cartier. 
In terms of bunches, we have the following characterization.

\begin{coro}
Suppose that $\Theta$ is free. Then $X$ 
is $\QQ$-Gorenstein if and only if we have
$$ 
\sum_{i=1}^{r} w_{i} 
\; \in \; 
\bigcap_{\tau \in \Theta} \lin(\tau).
\qquad \qed
$$
\end{coro}

Similarly we can decide whether a toric variety is a Fano
variety, i.e.~its anticanonical divisor is Cartier and 
ample:

\begin{coro}\label{fanocrit}
\begin{enumerate}
\item Suppose that $\Theta$ is free. 
Then $X$ is Fano if and only if we have 
$$ 
\sum_{i=1}^{r} w_{i} 
\; \in \; 
\bigcap_{\gamma_{0} \in \cov(\Theta)} Q(\gamma_{0} \cap E)^{\circ} .
$$
\item Suppose that $\Theta$ is free and $X$ is smooth.
Then $X$ is Fano if and only if we have 
$$ 
\sum_{i=1}^{r} w_{i} 
\; \in \; 
\bigcap_{\tau \in \Theta} \tau^{\circ}.
\qquad \qed
$$
\end{enumerate}
\end{coro}

We now apply our results to the case of simple bunches. 
Recall that these
describe precisely the $\QQ$-factorial full
2-complete toric varieties.
The first observation is that in this case we immediately
obtain the Picard group.

\begin{prop}
Suppose that $\Theta$ is a simple bunch. 
Then $\mathfrak{D} \colon E \to \WDiv^{T}(X)$ 
induces an isomorphism $K \cong \Pic(X)$. 
\end{prop}

\proof
It suffices to show that the image of $\mathfrak{D}$ 
is the group of invariant Cartier divisors $\CDiv^{T} (X)$. 
As $\Theta$ is simple, Lemma~\ref{surjectiveandinjective2} 
tells us that the dual projected cone
$(F \topto{P} N, \delta)$ and the maximal projectable fan
$\Sigma$ associated to $\Theta$
are the  universal reduced Cox 
construction of the fan $\Delta$ of $X$.

Now, for an element $\rq{w} \in E$, consider the restrictions
$\rq{w}|_{\lin(\rq{\sigma}) \cap F}$, where
$\sigma \in \Delta^{\max}$. By the
properties of a universal reduced Cox construction,
$$P \colon \lin(\rq{\sigma}) \cap F \to N$$ 
is an isomorphism for each $\sigma \in \Delta^{\max}$. 
Hence, $\rq{w}|_{\lin(\rq{\sigma}) \cap F}$
defines a $\ZZ$-valued linear function on $N$.
This shows that $\mathfrak{D}(\rq{w})$ is in fact a Cartier 
divisor.

On the other hand, if $D$ is a Cartier divisor, 
we have to show that $\rq{w} := \mathfrak{D}^{-1} (D)$ is
$\ZZ$-valued on $F$. 
As above, we see that the
restriction of $\rq{w}$ to any sublattice $\lin(\rq{\sigma}) \cap F$ is
$\ZZ$-valued. 
This implies the assertion, because 
Lemma~\ref{freeexists}~(ii) yields
$$ F = \sum_{\rq{\sigma}} \lin(\rq{\sigma}) \cap F. \qquad \qed$$

For a $\QQ$-factorial toric variety over $\CC$, we can identify 
$H^{2}(X,\QQ)$ with $K_{\QQ}$.  
Hence we may identify $H_{2}(X, \QQ)$ with 
$L_{\QQ} = \Hom(K_{\QQ},\QQ)$.
Then the Mori Cone, i.e., the cone ${\rm NE}(X) \subset L_{\QQ}$ 
of numerically effective curve classes is 
dual to the cone of numerically effective divisor classes
$\mathcal{N} (X) \subset K_{\QQ}$.

Reid proved that the Mori Cone of a $\QQ$-factorial complete
toric variety is generated by the classes
of the invariant curves, see~\cite[Cor.~1.7]{Re}. 
Using the fact that on a $\QQ$-factorial complete toric variety 
semiample and numerically effective divisors coincide,
Theorem~\ref{QCartier} gives a new description of the Mori Cone:
 
\begin{coro}\label{MoriCone}
Suppose that $X_{\Theta}$ is complete and simplicial.
Then the cone of numerically effective 
curve classes in $H_{2}(X,\QQ)$ is given by
$$ 
{\rm NE}(X_{\Theta}) \cong  \sum_{\tau \in \Theta} \tau^{\vee}.
$$
In particular, this cone is convex and polyhedral. Moreover,
$X_{\Theta}$ is projective if and only if ${\rm NE}(X_{\Theta})$
is strictly convex.
\endproof
\end{coro}

\section{Applications and examples}\label{section10}

In this section, we present some applications and examples.
First, we perform Kleinschmidt's classification in the setting of
bunches. We use the visualization techniques introduced 
at the end of Section~\ref{section2}.

\begin{prop}\label{Kleinschmidt}
The 2-complete smooth toric varieties $X$ with $\Cl(X) \cong \ZZ^{2}$ 
and $\mathcal{O}(X) \cong \KK$ correspond to free bunches 
$\Theta = \{\cone(w_{1},w_{2})\}$ given by
\begin{itemize}
\item weight vectors $w_{1} := (1,0)$ and
  $w_{i} := (b_{i},1)$, with $0 = b_{n} < b_{n-1} < \cdots < b_{2}$, 
\item multiplicities $\mu_{i} := \mu(w_{i})$ with $\mu_{1} > 1$,
  $\mu_{n} > 0$ and $\mu_{2} + \cdots + \mu_{n} > 1$. 
\end{itemize}
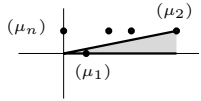
\begin{center}
  \begin{picture}(0,0)%
\includegraphics{smoothclass.pstex}%
\end{picture}%
\setlength{\unitlength}{1243sp}%
\begingroup\makeatletter\ifx\SetFigFont\undefined%
\gdef\SetFigFont#1#2#3#4#5{%
  \reset@font\fontsize{#1}{#2pt}%
  \fontfamily{#3}\fontseries{#4}\fontshape{#5}%
  \selectfont}%
\fi\endgroup%
\begin{picture}(3847,1937)(1126,-2333)
\put(2476,-1861){\makebox(0,0)[lb]{\smash{\SetFigFont{6}{7.2}{\familydefault}{\mddefault}{\updefault}
$(\mu_{1})$
}}}
\put(1126,-961){\makebox(0,0)[lb]{\smash{\SetFigFont{6}{7.2}{\familydefault}{\mddefault}{\updefault}
$(\mu_{n})$
}}}
\put(4101,-636){\makebox(0,0)[lb]{\smash{\SetFigFont{6}{7.2}{\familydefault}{\mddefault}{\updefault}
$(\mu_{2})$
}}}
\end{picture}

\end{center}
Moreover, the toric variety $X$ defined by such a bunch $\Theta$
is always projective, and it is Fano if and only if we have
$$ b_{2}(\mu_{3} + \cdots + \mu_{n}) < \mu_{1} + b_{2}\mu_{3} +
\cdots + b_{n-1}\mu_{n-1}. $$
\end{prop}

\proof 
We first show that
every smooth 2-complete toric variety $X$ with 
$\Cl(X) \cong \ZZ^{2}$ and $\mathcal{O}(X) \cong \KK$
arises from a bunch as in the assertion.
{From} Theorem~\ref{standard2maximal}~(ii), we infer
that the toric variety $X$ arises from a free bunch
$\Theta$ in a projected cone 
$(\ZZ^{m} \topto{Q} \ZZ^{2}, \gamma)$
with $\gamma = \cone(e_{1}, \ldots, e_{m})$.
Let $\{w_{1}, \ldots, w_{n}\}$ be the set of weight 
vectors, and let $\mu_{i}$ be the multiplicity of 
$w_{i}$.

Smoothness of $X$ means that every image
$Q(\gamma_{0} \cap \ZZ^{m})$, 
where $\gamma_{0} \in \cov(\Theta)$, 
generates $\ZZ^{2}$,
see Proposition~\ref{smoothchar}. 
Since we have $\mathcal{O}(X) = \KK$, 
Proposition~\ref{nontrivfunctionschar}
tells us that the cone $\vartheta \subset K_{\QQ}$
generated by the weight vectors is 
strictly convex. 
Being two-dimensional, $\vartheta$ 
is generated by two vectors, 
say $\vartheta = \cone(w_{1}, w_{n})$.

Now we make essential use of the fact that $\Theta$ 
lives in a two-dimensional space: Consider the intersection 
$\tau_{0}$ of all $\tau \in \Theta$.
Then $\tau_{0}$ is a weight cone. 
Moreover, $\tau_{0}$ is of dimension two,
because for any two cones of $\Theta$ their relative 
interiors intersect. 
Thus the defining property of a bunch implies 
$\Theta = \{\tau_{0}\}$.

As a strictly convex two-dimensional weight cone containing
$\tau_{0}$, also $\vartheta$ occurs among the images 
$Q(\gamma_{0})$, where $\gamma_{0} \in \cov(\Theta)$.
Consequently, $\vartheta = \cone(w_{1}, w_{n})$ is a regular 
cone in $\ZZ^{2}$,  
and moreover, using Remark~\ref{downstairsisochar}, we may
assume that $w_{1} = e_{1}$ and $w_{n} = e_{2}$ are 
the canonical base vectors of
$\ZZ^{2}$.

We claim that the cone $\tau_{0}$ has at least one ray 
in common with $\vartheta$.
Indeed, otherwise one of the cones generated
by $\tau_{0} \cup \QQ_{\ge 0} e_{1}$ or 
$\tau_{0} \cup \QQ_{\ge 0} e_{2}$ would be 
non regular.
But both cones occur among $Q(\gamma_{0})$, where 
$\gamma_{0} \in \cov(\Theta)$. A contradiction.
Again by Remark~\ref{downstairsisochar}, we may
assume that $\tau_{0}$ contains $e_{1}$, 
and hence is generated by vectors $e_{1}$
and $b_{2}e_{1} + e_{2}$, where $b_{2} \ge 0$.

Consider any $w_{i} \ne e_{1}$. Then 
$\tau_{i} := \cone(w_{i},e_{1})$ overlaps $\tau_{0}$,
and thus the properties of a bunch give 
$\tau_{0} \subset \tau_{i}$.
Consequently $\tau_{i} = Q(\gamma_{0})$ 
for some $\gamma_{0} \in \cov(\Theta)$,
and hence $w_{i}$ and $e_1$ generate $\ZZ^{2}$. 
But this means $w_{i} = b_{i}e_{1} + e_{2}$
with $b_{i} \ge b_{2}$.  
So we arrive at the desired picture.  
Note that the conditions $\mu_{1} > 1$ 
and $\mu_{2} + \ldots + \mu_{n} > 1$ are due
to Property~\ref{standarddefi}~(ii).

The fact that the toric variety $X$ associated to $\Theta$ is
projective is clear by Proposition~\ref{completenesschar}
and Corollary~\ref{Ewaldcrit}. 
Moreover, Corollary~\ref{fanocrit}
tells us that $X$ is Fano if and only if the weighted sum
$\mu_{1}w_{1} + \ldots + \mu_{n}w_{n}$ lies in the relative
interior of the cone $\tau$. But this holds if and only if the
condition stated in the assertion is valid. \endproof

Let us illustrate the quasiprojectivity criterion~\ref{Ewaldcrit} 
by means of the following version of a classical example taken from 
Oda's book~\cite[p.~84]{Od}:

\begin{exam}\label{oda}
The simplest nonprojective complete simplicial fans $\Delta$ live in
$\ZZ^{3}$ and are of the following type: 
Consider a prism $P \subset \QQ^{3}$ over a 2-simplex
such that $0$ lies in the relative interior of $P$. 
Here we take $P$ with the vertices
$$
(-1,0,0), \quad
(0,-1,0), \quad
(0,0,-1), \quad
(0,1,1), \quad 
(1,0,1), \quad
(1,1,0). 
$$
Then subdivide the facets of $P$ according to the picture below, 
and let $\Delta$ be the fan generated by the cones over the simplices
of this subdivision.
\begin{center}
  \begin{picture}(0,0)%
\includegraphics{oda_cut1.pstex}%
\end{picture}%
\setlength{\unitlength}{2072sp}%
\begingroup\makeatletter\ifx\SetFigFont\undefined%
\gdef\SetFigFont#1#2#3#4#5{%
  \reset@font\fontsize{#1}{#2pt}%
  \fontfamily{#3}\fontseries{#4}\fontshape{#5}%
  \selectfont}%
\fi\endgroup%
\begin{picture}(9135,3406)(676,-2612)
\put(8371,569){\makebox(0,0)[lb]{\smash{\SetFigFont{8}{9.6}{\rmdefault}{\mddefault}{\updefault}
$\frac{1}{2}w_3$
}}}
\put(7201,389){\makebox(0,0)[lb]{\smash{\SetFigFont{8}{9.6}{\rmdefault}{\mddefault}{\updefault}
$e_3$
}}}
\put(9811,209){\makebox(0,0)[lb]{\smash{\SetFigFont{8}{9.6}{\rmdefault}{\mddefault}{\updefault}
$e_2$
}}}
\put(9361,-421){\makebox(0,0)[lb]{\smash{\SetFigFont{8}{9.6}{\rmdefault}{\mddefault}{\updefault}
$\frac{1}{2}w_1$
}}}
\put(9001,-961){\makebox(0,0)[lb]{\smash{\SetFigFont{8}{9.6}{\rmdefault}{\mddefault}{\updefault}
$e_1$
}}}
\put(7201,-286){\makebox(0,0)[lb]{\smash{\SetFigFont{8}{9.6}{\rmdefault}{\mddefault}{\updefault}
$\frac{1}{2}w_2$
}}}
\put(676,-2536){\makebox(0,0)[lb]{\smash{\SetFigFont{9}{10.8}{\familydefault}{\mddefault}{\updefault}
Defining polytope subdivision of $\Delta$
}}}
\put(6976,-2536){\makebox(0,0)[lb]{\smash{\SetFigFont{9}{10.8}{\familydefault}{\mddefault}{\updefault}
Corresponding bunch $\Theta$
}}}
\end{picture}

\end{center}
The bunch $\Theta$ corresponding to $\Delta$ is free and lives 
in a 3-dimensional lattice $K \cong \ZZ^{3}$. 
Combinatorially, it looks as indicated above. 
Explicitly it is given by the weight vectors
$e_{1}$, $e_{2}$, $e_{3}$, 
$w_{1} := e_{1} + e_{2}$, 
$w_{2} := e_{1} + e_{3}$, 
and $w_{3} = e_{2} + e_{3}$ in $\ZZ^{3}$, 
and the four cones
$$ 
\begin{array}{ll}
\cone(e_{3}, w_{1}, w_{2}), \quad & \cone(e_{1}, w_{1}, w_{3}) \\
\cone(e_{2}, w_{2}, w_{3}), \quad & \cone(w_{1}, w_{2}, w_{3})
\end{array}$$

Using the Propositions~\ref{fullandsimplicial},
\ref{completenesschar}, and Corollary~\ref{Ewaldcrit} 
it is immediately clear from the picture that 
$X := X_{\Theta}$ is $\QQ$-factorial and nonprojective.
Moreover, by Theorem~\ref{QCartier}, the cone $\C^{\rm sa} (X)$ 
of semiample divisors is spanned by the class of the 
anticanonical divisor.
\end{exam}

The last example concerns the problem whether or not
the Betti numbers of a toric variety are determined by
the combinatorial type of the defining fan.
A first counterexample was given by McConnell, see~\cite{MC}.
Later Eikelberg~\cite{Ei} gave the following simpler one:

\begin{exam}\label{eikelberg}
Let $P \subset \QQ^{3}$ be a prism over a 2-simplex such that
$0$ lies in the relative interior of $P$,
and define $\Delta$ to be the fan generated the cones over the 
facets of $P$. As vertices of $P$ we take:
$$
\begin{array}{lll}
v_{1} := (1,0,1), &
v_{2} := (0,1,1), &
v_{3} := (-1,-1,1), \\
v_{4} := (1,0,-1), &
v_{5} := (0,1,-1), &
v_{6} := (-1,-1,-1).
\end{array}
$$
Eikelberg constructs a nonprojective fan $\Delta'$ from $\Delta$ by 
moving the ray through $v_{2}$ into the ray
through $v_{2}' := (1,2,3)$.
The picture is the following 
(the dotted diagonals indicate the edges 
of the convex hull over $v_{1},v_{2}',v_{3},\ldots,v_{6}$
that do not define cones of $\Delta'$):
\begin{center}
  \begin{picture}(0,0)%
\includegraphics{betti_fan.pstex}%
\end{picture}%
\setlength{\unitlength}{2072sp}%
\begingroup\makeatletter\ifx\SetFigFont\undefined%
\gdef\SetFigFont#1#2#3#4#5{%
  \reset@font\fontsize{#1}{#2pt}%
  \fontfamily{#3}\fontseries{#4}\fontshape{#5}%
  \selectfont}%
\fi\endgroup%
\begin{picture}(9044,3473)(1329,-2837)
\put(7651,-2761){\makebox(0,0)[lb]{\smash{\SetFigFont{9}{10.8}{\familydefault}{\mddefault}{\updefault}
One ray moved
}}}
\put(1801,-2761){\makebox(0,0)[lb]{\smash{\SetFigFont{9}{10.8}{\familydefault}{\mddefault}{\updefault}
Ordinary prism $P$
}}}
\end{picture}

\end{center}

What does this mean in terms of bunches? First note that the
bunches $\Theta$ and $\Theta'$ corresponding to $\Delta$ and $\Delta'$
contain 2-dimensional cones, because $\Delta$ as well as $\Delta'$
have nonsimplicial  cones.
The loss of projectivity is reflected in terms of bunches as follows:
\begin{center}
  \begin{picture}(0,0)%
\includegraphics{betti_bunch.pstex}%
\end{picture}%
\setlength{\unitlength}{2072sp}%
\begingroup\makeatletter\ifx\SetFigFont\undefined%
\gdef\SetFigFont#1#2#3#4#5{%
  \reset@font\fontsize{#1}{#2pt}%
  \fontfamily{#3}\fontseries{#4}\fontshape{#5}%
  \selectfont}%
\fi\endgroup%
\begin{picture}(9530,3429)(1126,-2961)
\put(1126,-2885){\makebox(0,0)[lb]{\smash{\SetFigFont{9}{10.8}{\familydefault}{\mddefault}{\updefault}
Bunch $\Theta$ corresponding to $\Delta$
}}}
\put(6526,-2885){\makebox(0,0)[lb]{\smash{\SetFigFont{9}{10.8}{\familydefault}{\mddefault}{\updefault}
Bunch $\Theta'$ corresponding to $\Delta'$
}}}
\end{picture}

\end{center}

Here we draw the intersection of $\Theta$ and $\Theta'$ with a plane
orthogonal to an inner vector of the (strictly) convex hulls
$\vert \Theta \vert$ and $\vert \Theta' \vert$ of the unions of the
respective bunch cones. 
So the fat lines correspond to 2-dimensional bunch cones, 
whereas the shaded simplices represent 3-dimensional bunch cones.

Now, both fans $\Delta$ and $\Delta'$ have the same combinatorial
type. In terms of bunches, we see immediately that the associated 
toric varieties $X$ and $X'$ have different Betti numbers:
The second Betti number $b_{2}(X)$ equals the dimension of 
$\Pic_{\QQ}(X)$.
Hence Theorem~\ref{QCartier} gives us $b_{2}(X) = 1$.
For $X'$, we obtain $b_{2}(X') = 0$ by the same reasoning.
\end{exam}

We conclude with a remark concerning
Ewald's construction of ``canonical extensions''
of a given toric variety presented in~\cite[Section~3]{Ew},
compare also~\cite{Bon}.
The bunch theoretical analogue is the following:

\begin{constr}
Consider a free bunch, represented in the sense of
Construction~\ref{bunchfromdownstairs}
by a set weight vectors $\{w_{1}, \ldots, w_{n}\}$ 
in a lattice $K$, 
multiplicities $\mu_{1}, \ldots, \mu_{n}$ 
and a collection $\Theta$ of weight cones.
Define a new free bunch by setting 
$$
K' := K, \qquad
w_{i}' := w_{i}, \qquad
\Theta' := \Theta
$$
and replacing the multiplicities $\mu_{i}$ with bigger
ones, say $\mu_{i}'$. 
For the toric varieties $X$ and $X'$ associated to these 
bunches, we have 
$$
\dim(X') - \dim(X) = \sum_{i=1}^{n} (\mu_{i}' - \mu_{i}),
\qquad
\ClDiv(X') \cong \ClDiv(X),
\qquad
\Pic(X') \cong \Pic(X).
$$ 
Moreover, applying the respective characterizations
of Section~\ref{section7} and \ref{section9}, 
one immediately verifies that $X'$ is non quasiprojective 
(complete, $\QQ$-factorial, smooth) if $X$ was so.
\end{constr}

\end{document}